\documentclass[11pt]{article}
\usepackage{amssymb}

\usepackage{color}

\usepackage{epsfig}
\newtheorem{Lemma}{Lemma}[section]
\newtheorem{Theorem}[Lemma]{Theorem}
\newtheorem{Proposition}[Lemma]{Proposition}
\newtheorem{Corollary}[Lemma]{Corollary}
\newtheorem{Remark}[Lemma]{Remark}
\newtheorem{Remarks}[Lemma]{Remarks}
\newtheorem{Definition}[Lemma]{Definition}
\newtheorem{Assumption}[Lemma]{Assumption}
\newtheorem{Example}[Lemma]{Example}

\def\theequation{\arabic{section}.\arabic{equation}}
\def\nc{\newcommand}
\nc{\mylabel}[1]{\label{#1}}
\nc{\isep}{\itemsep -0.6mm}
\nc{\ds}{\displaystyle}

\nc{\rr}    {\rightarrow}

  \nc{\tf}{{\tilde{f}}}
  \nc{\tg}{{\tilde{g}}}
\nc{\tx}{{\tilde{x}}} \nc{\ty}{{\tilde{y}}}
\nc{\tF}{{\tilde{F}}}

\nc{\us}{\underline{s}}    \nc{\os}{\overline{s}}
\nc{\UP}{\underline{P}}    \nc{\OP}{\overline{P}}
\nc{\OD}{\overline{D}}     \nc{\OB}{\overline{B}}
\nc{\HD}{\mbox{${\rm dim}_H$}}
\nc{\dimH}{{\rm dim}_H}

\renewcommand{\Re}    {{\rm Re}}
\renewcommand{\Im}    {{\rm Im}}

\newcommand{\ess}    {{\rm ess}}
\newcommand{\osc}    {{\rm osc}}

\nc{\hra}{\hookrightarrow}
 \nc{\Arg}   {{\rm Arg}}

\nc{\grass}{{{\rm Gr}_2}} \nc{\dist}{{\rm dist}}
\nc{\Lip} {{\rm Lip}}
\nc{\Card}{{\rm Card}}
\nc{\Conv} {{\rm Conv}}
\nc{\Span}{{\rm Span}}
\nc{\SpanC}{{\rm Span}_\CC}
\nc{\SpanR}{{\rm Span}_\RR}
\nc{\re}{{\rm Re}} \nc{\im}{{\rm Im}} \nc{\Aut}{{\rm Aut}}
\nc{\Area}{{\rm Area}} \nc{\size}{{\rm size}} \nc{\crit}{{\rm crit}}
\nc{\const}{{\rm const}} \nc{\diam}{{\rm diam}} \nc{\bfh}{{\bf h}}
\nc{\bfone}{{\bf 1}}  

 \nc{\Int}{{\rm Int\;}} \nc{\Imm}{\mbox{Im}}
\nc{\half}{\frac{1}{2}} \nc{\DDD}{D}

\nc{\calP}{{\cal P}}
\nc{\calA}{{\cal A}}
\nc{\calC}{{\cal C}}
\nc{\CS}{{\cal C}^*}
\nc{\ellone}{{\{\ell=1\}}}
\nc{\Cellone}{{{\cal C} \cap \ellone}}
\nc{\D}{{\cal D}} \nc{\EE}{{\cal E}}

\nc{\CCR}{{{\cal C}_{\RR}}} 
\nc{\CCSR}{{{\cal C}_{\RR}^*}} 
\nc{\CCPR}{{{\cal C}^{'}_{\RR}}} 
\nc{\CCOR}{{{\cal C}^1_{\RR}}} 
\nc{\dCCR}{d_\CCR}
\nc{\dCC}{d_{\cal C}}
\nc{\CCS}{{{\cal C}^*}} 
\nc{\CCSC}{{{\cal C}_{\CC}^*}} 
\nc{\CCC}{{{\cal C}_{\CC}}} 
\nc{\CCOC}{{{\cal C}^1_{\CC}}} 
\nc{\dCCC}{d_\CCC}
\nc{\CCone}  {{{\CC}^n_+}}
\nc{\XP}   {X^{'}}

\nc{\CC}{{\Bbb C}} 
\nc{\CCbf}{{\Bbb C}} 
\nc{\DD}{{\Bbb D}} \nc{\EEE}{\mbox{$\Bbb E$}}
\nc{\HH}{{\Bbb H}}
\nc{\RR}{{{\Bbb R}}}
\nc{\RRs}{{\Bbb R}} \nc{\NN}{\mbox{$\Bbb N$}} \nc{\NNs}{{\Bbb N}}
\nc{\ZZ}{\mbox{$\Bbb Z$}} \nc{\AutDD}    {{\rm Aut(\DD;\diag \DD)}}
\nc{\hatRR}{\widehat{\RR}}
\nc{\hatCC}{\widehat{\CC}}

\nc{\la}{\langle}
\nc{\ra}{\rangle}
\nc{\barA}{\overline{A}} 
\nc{\barM}{\overline{M}} 
\nc{\bara}{\overline{a}} 
\nc{\barz}{\overline{z}} 
\nc{\baru}{\overline{u}} 
\nc{\barb}{\overline{b}} \nc{\barDD}{{\overline{\DD}}}
\nc{\barK}{{\overline{K}}} \nc{\bRs}{\overline{\RR}_+((s))}
\nc{\UbarU}{U\times \overline{U}} \nc{\XbarX}{X\times \overline{X}}
\nc{\barmu}{\overline{\mu}}
\nc{\barla}{\overline{\lambda}}

\nc{\hatD}{\widehat{D}}
\nc{\hatx}{\widehat{x}}
\nc{\haty}{\widehat{y}}
\nc{\hatd}{\widehat{d}} \nc{\hatg}{g_\htinyD} \nc{\hatv}{\widehat{v}}
\nc{\Dhatf}{D\widehat{f}} \nc{\Dpsihat}{D\widehat{\psi}_t}
\nc{\hatf}{{\widehat{f}}}

\nc{\Cl}{{\rm Cl\;}} \nc{\len}{\mbox{len\,}} \nc{\diag}{\mbox{diag\ }}
\nc{\diagK}{\mbox{diag}(K)}

\nc{\Halmos}    {\ \raisebox{0.3ex}  {\framebox[0.9ex]{
			   \rule[0ex]{0ex}{0.5ex}
			    }}}
\nc{\longr}{\longrightarrow}

\begin{document}
\title{ Cones and gauges in complex spaces~:\\
 Spectral gaps  and complex Perron-Frobenius theory.
%\footnote{cc-81}
}
\author{Hans Henrik Rugh\\
    University of Cergy-Pontoise, France}
\date {\today}
 \maketitle
\begin{abstract}
We introduce complex cones and associated projective gauges,
  generalizing a real Birkhoff cone and its
Hilbert metric to complex vector spaces.
 We deduce a variety of
{\em spectral gap} theorems in complex Banach spaces.
We prove a {\em dominated} complex cone-contraction Theorem
 and use it to extend the classical
  Perron-Frobenius Theorem to complex matrices,
 Jentzsch's Theorem to complex integral operators, 
a Kre\u{\i}n-Rutman Theorem to compact and quasi-compact complex operators
and a Ruelle-Perron-Frobenius Theorem to complex transfer operators
in dynamical systems.
In the simplest case of a complex $n$ by $n$ matrix
$A\in M_n(\CC)$ we have the following statement~:
Suppose that 
 $0<c<+\infty$ is such that
 $ |\im\, A_{ij}\overline{A}_{mn}|
    <  c \leq \re\, A_{ij}\overline{A}_{mn}$ for
all indices.  Then $A$ has a `spectral gap'. 
\end{abstract}

\section{Introduction}
 The Perron-Frobenius Theorem \cite{Per07,Fro08}
 asserts that
 a real square matrix with
 strictly positive entries has a `spectral gap', i.e.\
 the matrix has a positive simple eigenvalue and 
 all other eigenvalues are strictly smaller in modulus.
 More generally, let $A$ be a bounded linear operator acting
upon a real or complex Banach space and of spectral radius
$r_{\rm sp}(A)$.
We say that $A$ has  a spectral gap if (1)
 it has  a simple isolated
eigenvalue $\lambda$ the modulus of which equals 
$r_{\rm sp}(A)$
 and (2) the remaining part  of the spectrum is contained
in a disk centered at zero and of radius strictly smaller than 
$r_{\rm sp}(A)$.

Jentzsch generalized  in \cite{Jen12}
the Perron-Frobenius Theorem to integral operators with a strictly positive
continuous kernel. 
The proof uses the Schauder-Tychonoff Theorem
 to produce a dual eigenvector
and then a contraction on the kernel of this eigenvector
 to get a spectral gap.
Kre\u{\i}n-Rutman
\cite[Theorem 6.3]{KR50} (see also \cite{Rut40} and 
   \cite{Rot44}) gave
   an abstract setting of this result by considering
a punctured real closed cone mapped to its interior by a compact operator.
Compactness of the operator essentially reduces the problem to
finite dimensions.

 Birkhoff,
in a seminal
paper \cite{Bir57},
% (and also \cite{Bir67}),
 developed a more elementary and intuitive (at least in our opinion)
 Perron-Frobenius `theory'
by considering the projective contraction of a cone
 equipped with its associated  Hilbert metric.
Birkhoff noted that
this projective metric satisfies a  {\em contraction principle},
 i.e.\ any linear map
preserving the cone is a contraction for the metric and the contraction
is  strict and uniform if the image
of the cone has finite projective diameter.

All these  results, or rather their proofs,
make use of the `lattice'-structure
induced by a real cone on a real Banach space 
(see \cite{Bir67} and also \cite{Mey91}).
 On the other hand, from complex analysis we know that the 
Poincar\'e metric on the unit disk, 
$\DD=\{z\in \CC: |z|<1\}$, and the induced metric
on a hyperbolic Riemann surface enjoy properties similar
to the Hilbert metric, in particular a {\em contraction principle}
with respect to conformal maps.  More precisely,
 if $\phi: U\rr V$ is a conformal map between
 hyperbolic Riemann surfaces 
 then its conformal derivative never exceeds one.
The map is a strict contraction unless it
  is a bijection (see  e.g.\ 
\cite[Chapter I.4: Theorems 4.1 and 4.2]{CG93}).
By considering analytic images of complex discs
 Kobayashi \cite{Kob67,Kob70} (see also \cite{Ves76})
 constructed  a hyperbolic  metric on 
complex (hyperbolic) manifolds, a tool with many applications
also in  infinite dimensions (see e.g.\ \cite[Appendix D]{Rug02}).
%Vesentini \cite[section 6]{Ves76} `complexified' a real cone
 %$\CCR\subset X$ by considering a `tube' domain
 %$X \oplus i\CCR$.  He used this to construct a 
%Kobayashi like metric on the direct sum. His complexification is, however,
%not $\CC$-invariant and this prevented him from going beyond 
%real operator contractions.

Given a real cone contraction,
perturbation theory allows on abstract grounds to consider
`small' complex perturbations but uniform estimates are usually
 hard to obtain.  Uniform complex estimates are needed e.g.\ when
proving local limit theorems and refined large deviation theorems
for Markov additive processes (see \cite{NN87} and references therein)
and also for studying the regularity of
characteristic exponents for time-dependent
and/or random dynamical systems (see e.g.\ \cite{Rue79,Rug02}).
 It is  desirable to obtain a description of a
projective contraction and, in particular, a spectral gap condition
 for complex operators without the above-mentioned  restrictions. 
We describe in the following one way to accomplish this goal.\\
 
In section \ref{sec C cones}
we introduce families of $\CC$-invariant cones
in complex Banach spaces and a theory for the projective contraction
of such cones.
   %(see Figure \ref{fig cone ex}, section \ref{sec complex} for 
   %an illustration).
 %abandon convexity as well as
  %the partial ordering that a
%vector space otherwise inherits
 %This, however, does not affect
%the essential properties needed for a contraction.
 The central idea is simple, namely 
to use the Poincar\'e metric as a `gauge' on 2-dimensional
affine sections of a complex cone. 
At first sight, this looks like the Kobayashi construction.
A crucial difference, however, is that we only consider disk images 
in 2-dimensional subspaces. Also we do not take infimum over chains
(so as to obtain a triangular inequality, see
  Appendix \ref{a proj metric}).
 This adapts well to the study
of  linear operators and makes computations 
much easier than for the general Kobayashi metric.
Lemma \ref{projective gauge} shows that
this gauge is indeed projective.  The contraction
principle for the Poincar\'e metric  translates into a contraction
principle for the gauge and, under additional regularity assumptions,
developed in section \ref{sec complex}, into a projective contraction,
and finally a spectral gap, with respect to the Banach space norm.

In sections 
\ref{section real case} and \ref{sec canonical complexification}
we consider real cones and define
their {\em canonical complexification}.  For example,
$\CC_+^n = \{ u\in \CC^n : |u_i+u_j|\geq |u_i-u_j|, \forall\, i,j\}=
 \{ u\in \CC^n : \Re\, u_i\baru_j\geq 0, \forall\, i,j\} $
is the canonical complexification of the standard
real cone, 
$\RR_+^n$.
We show that our complex cone contraction yields  a genuine extension of the
Birkhoff cone contraction~: A real Birkhoff cone is 
{\em isometrically}
embedded into its canonical  complexification.
 It enjoys here the same contraction properties with respect to linear
 operators. We obtain in section \ref{C dominated}
then one of our main results:
When a complex operator is
{\em dominated} by a sufficiently regular real cone-contraction
(Assumption \ref{main assumption})
then (Theorem
\ref{C domin sp gap}) the complex operator has a
spectral gap. It is of interest to note
 that the conditions on the complex operator
are expressed in terms of a real cone and sometimes
easy to verify.
Sections \ref{sec Applications}-\ref{sec RPF}
	thus presents a selection of complex 
analogues of well-known real cone contraction theorems~:
 A Perron-Frobenius Theorem
for complex matrices (as stated at the end of the abstract),
 Jentzsch's Theorem for  complex integral operators,
 a Kre\u{\i}n-Rutman Theorem for compact and quasi-compact
 complex operators and
a Ruelle-Perron-Frobenius Theorem for complex
transfer operators. 

Acknowledgements: I am  grateful to A Douady for a suggestion
in the proof of
Lemma \ref{lemma contraction} and to the anonymous referee for several
valuable suggestions and corrections.

\section{Complex cones and gauges}
 \mylabel{sec C cones}
Let $\hatCC=\CC\cup \{\infty\}$ denote 
the Riemann sphere. When  $U\subset\hatCC$ is an open connected
 subset avoiding at least three points one says that
 the set is {\em hyperbolic}.
We write $d_U$ for the corresponding 
hyperbolic metric.
We refer to \cite[Chapter I.4]{CG93} or 
\cite[Chapter 2]{Mil99} for 
the properties of the hyperbolic metric which we use in the present paper.
 As normalization we use 
 $ds = 2|dz|/(1-|z|^2)$ on the unit disk $\DD$
 and the metric $d_U$ on $U$ induced
by a Riemann mapping $\phi: \DD \rr U$. One then has~:
\begin{equation}
         d_{\DD}(0,z) =  \log \frac{1+|z|}{1-|z|}, \ \ \ \
        |z| = \tanh \frac{d_\DD(0,z)}{2}.
    \mylabel{C metric}
 \end{equation}

Let  $E$ be a  complex topological vector space.
We denote by $\Span\{x,y\}=\{\lambda x+ \mu y : \lambda,\mu\in\CC\}$
the complex subspace generated by two vectors $x$ and $y$ in $E$.

\begin{Definition}
\mylabel{def complex cone}

\mbox{}
\begin{enumerate}
\item[(1)] We say that a subset $\calC\subset E$
is a {\bf closed complex cone} if it is closed in $E$,
$\CC$-invariant (i.e.\  $\calC=\CC \;\calC$) and
contains at least one complex line.
\item [(2)]
 We say that the closed complex cone $\calC$ is {\bf proper} if it
 contains no complex planes, i.e.\ if $x$ and $y$ are independent vectors
 then $\Span\{x, y\} \not\subset \calC$.
\end{enumerate}

Throughout this paper we will simply refer to a proper closed complex
cone as a {\bf $\mathbf \CCbf$-cone}. 
\end{Definition}

\noindent 
 Let  $\calC$ be a $\CC$-cone.
Given  a pair of non-zero vectors, $x,y\in \calC^*\equiv \calC-\{0\}$,
we consider the subcone~:
 $\Span\{x,y\}\cap \calC$. 
We wish to construct a `projective distance' between
the complex lines $\CC x$ and $\CC y$ within this subcone.
%(thus a distance in $\CC P^1$). 
We do this by considering
 the affine plane through $2x$ and $2y$, choosing
coordinates (the choice to some extend being arbitrary) as follows~: 
\begin{equation}
   D(x,y) \equiv D(x,y;\calC) = \{ \lambda\in\hatCC :
        (1+\lambda) x + (1-\lambda) y \in \calC\}\subset \hatCC,
	\mylabel{slice}
\end{equation}
with the convention that $\infty\in D(x,y)$ iff $x-y\in\calC$.
The interior of this "slice" is denoted  
$D^o(x,y)$ (for the spherical topology on $\hatCC$).
 We note that when $x$ and $y$ are linearly independent, continuity of
the canonical mapping
$\CC^2\rr \Span\{x,y\}$ implies that $D=D(x,y)$ is a closed subset
of $\hatCC$. As the cone is proper, $D\subset \hatCC$
 is a strict subset
so that $\hatCC-D$ is open and non-empty, whence contains (more than) 3 points.
If, in addition,
  $D^o$
is connected it is a hyperbolic
Riemann surface
 (\cite[Theorem I.3.1]{CG93}).

\begin{Definition}
\mylabel{def slicable cone}
%A proper complex cone, $\calC$, is said to be a
 %{\em  hyperbolic complex cone} if for every linearly independent vectors,
 %$x,y\in\calC^*$ the slice
    %$D(x,y)\subset \hatCC$ is hyperbolic (omits at least three points).
%
%We say that a hyperbolic cone is  {\em strongly hyperbolic}  if
%$D^o(x,y)$ is connected (when non-empty) for every $x,y\in\calC^*$.
%
Given a $\CC$-cone, we define the
   {\bf gauge}, $\dCC: \CS \times \CS \rr [0,+\infty]$,
 between two points
$x,y\in\calC^*$ as follows~:
When two vectors are co-linear 
%or, more generally, if $U$ is not hyperbolic
we set $\dCC(x,y)=0$.
If they are linearly independent and $-1$ and $1$ belongs to the same
 connected component $U$ of $D^o(x,y)$ 
 we set~:
   \begin{equation}
      d_\calC(x,y) \equiv d_{U}(-1,1) >0.
      \mylabel{gauge}
       \end{equation}
In all remaining cases,  we set 
 $\dCC=\infty$.

 When $V\subset \calC$ is a (sub-)cone of the $\CC$-cone $\calC$ 
 we write $\diam_\calC (V^*) \equiv \sup_{x,y\in V^*} \dCC(x,y)\in [0,+\infty]$
  for the
projective {\bf `diameter'} of $V$ in $\calC$.
We call it a diameter even though the gauge
 need not verify  the triangular inequality, whence
 need not be a metric
  (see Appendix
  \ref{a proj metric} for more on this issue).
\end{Definition}

\begin{Lemma}
\mylabel{projective gauge}
 Let $\calC$ be a  $\CC$-cone.
  The gauge on the cone is symmetric and projective, i.e.\ for 
 $x,y\in\calC^*$
  and $a\in\CC^*$~:
  \[ \dCC(y,x)=\dCC(x,y)=\dCC(ax,y)=\dCC(x,ay) .\]
\end{Lemma}
Proof:
For $(1+\mu) a + (1 - \mu) \neq 0$  we write 
  \[ (1+\mu) a x + (1-\mu) y = 
       \frac{(1+\mu) a + (1 - \mu)}{2} ((1+R)x + (1-R) y)\]
with
 \[ R = R_a(\mu) = \frac{ (1+\mu) a - (1-\mu) } 
        {(1+\mu) a + (1-\mu)} .\]
Then  $R_a$ extends to
a conformal bijection $R_a : \mu\in  D(ax,y)\mapsto R_a(\mu) \in D(x,y)$ 
(a M\"obius transformation of $\hatCC$)
 preserving $-1$ and $1$. The hyperbolic metric
is invariant under such transformations so
indeed $\dCC(x,y)=\dCC(ax,y)$ (but both could be infinite).
 Similarly, the map 
$\lambda\mapsto -\lambda$ yields a conformal bijection between 
the domains $D(x,y)$ and $D(y,x)$, interchanging $-1$ and $1$
and   the symmetry follows.
\Halmos\\

\begin{figure}
\begin{center}
\epsfig{figure=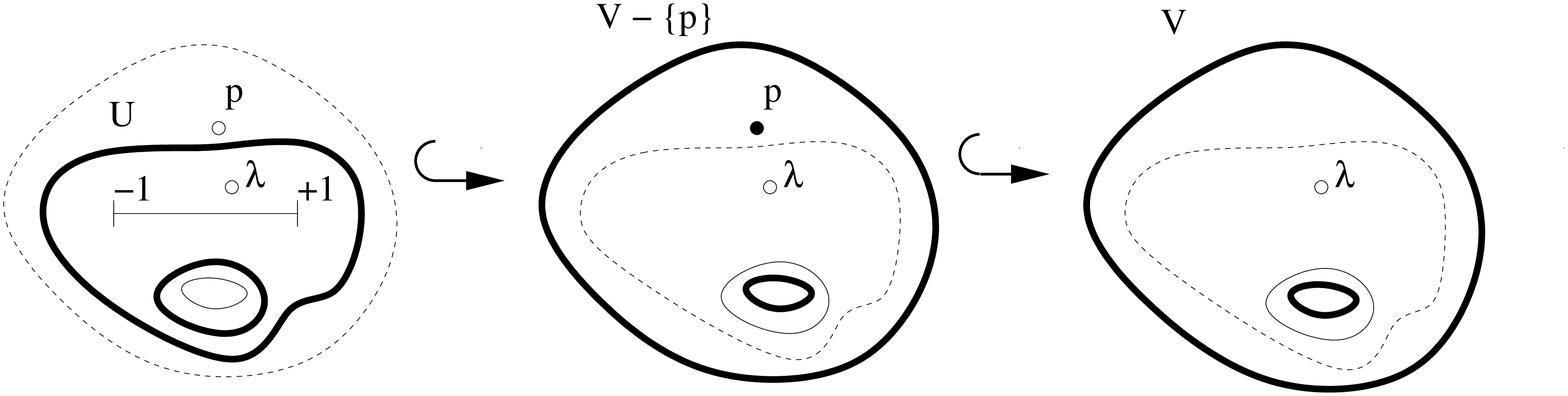,width=15cm}
\end{center}
\caption{The sequence of inclusions $U \hra V- \{p\} \hra V$ 
 in the proof of Lemma \ref{lemma contraction}}
\mylabel{fig inclus}
\end{figure}

\begin{Lemma}
\mylabel{lemma contraction}
Let $T:E_1\rr E_2$ be a complex linear map between topological  vector
spaces and let $\,\calC_1\subset E_1$ and
 $\,\calC_2\subset E_2$ be $\CC$-cones for which 
$T(\CS_1)\subset \CS_2$. Then the map, 
\[  T :
 (\CS_1,d_{\calC_1})
        \rr
 (\CS_2,d_{\calC_2}),
  \]
 is a contraction. 
If the image has finite diameter, i.e.\ 
$\Delta=\diam_{\calC_2^*} T \CS_1 < \infty$,  then the contraction is 
strict and uniform. More precisely, there is 
$\eta=\eta(\Delta)<1$ (depending on $\Delta$ only)
for which  
   \[ d_{\calC_2}(Tx,Ty)\leq \eta \ d_{\calC_1}(x,y), \ \ \ 
       \forall x,y\in\calC_1^* .\]
\end{Lemma}
Proof :
 Let $x,y\in\calC_1^*$ and set $D_1=D(x,y;\calC_1)$ 
 and $D_2=D(Tx,Ty;\calC_2)$ for which we have
 \[ \{-1,1\} \subset D_1 \subset D_2 \subset \hatCC. \]
  Suppose that $Tx,Ty\in\calC^*_2$ are linearly
independent and that $D_2$ and $D_1$ are  hyperbolic
 (if not, $d_{\calC_2}(Tx,Ty)$ vanishes and we are through).
Since shrinking a domain increases hyperbolic distances,
it follows that 
 $d_{\calC_2}(Tx,Ty)\leq  d_{\calC_1}(x,y)$ (although both could be infinite).

Suppose now that $\Delta<+\infty$. 
Then $-1$ and $1$ belong to the same
connected component, $V$, of $D^o(Tx,Ty)$. We may suppose that
 $-1$ and $1$ also belong to the same connected component, $U$,
of $D^o(x,y)$ (or else $d_{\calC_1}(x,y)=\infty)$ and we are through).
 Our assumptions imply that 
$U\subset V$ is a strict inclusion and that $\diam_V(U)\leq \Delta$.
 Choose 
 $\lambda\in U$ and
\cite{Dou04} pick
  $p\in V\setminus U$ for which
 $d_V(\lambda,p)\leq \Delta$
(this is possible as the inclusion $U\subset V$ is strict and the diameter
of $U$ did not exceed $\Delta$).
The inclusion $U\hookrightarrow V-\{p\}$ is non-expanding and
the inclusion $V-\{p\}\hookrightarrow V$
is a contraction which has conformal derivative uniformly smaller than
some $\eta=\eta(\Delta)<1$ on the  punctured  $\Delta$-neighborhood,
$B_V(p,\Delta)^*$, of p  (see Remark \ref{remark punctured disk}).
In particular, the composed map (see Figure \ref{fig inclus})
$U \hookrightarrow V-\{p\}\hookrightarrow V$ has conformal derivative
 smaller than $\eta(\Delta)$ at $\lambda\in B_V(p,\Delta)^*$.
 As $\lambda\in U$ 
was arbitrary this is true at any point along a
 geodesic joining $-1$ and $1$ in $U$ so that 
\[
   d_{\calC_2}(Tx,Ty)=
  d_{V}(-1,1) \leq 
   \eta \; d_{U}(-1,1)=
   \eta  \;d_{\calC_1}(x,y). \Halmos \] 
\mbox{}
\begin{Remark}
\label{remark punctured disk}
An explicit bound may be given using the expression
 $ds={|dz|}/{(|z| \log \frac{1}{|z|})}$ for the metric on the punctured 
 disk at $z\in\DD^*$
(see e.g.\ \cite[Example 2.8]{Mil99}).
Denoting, $t=\tanh \Delta/2$, 
 we obtain the bound,
  $\eta(\Delta)= \frac{2t}{1-t^2} \log \frac{1}{t}=
 \sinh (\Delta) \;
            \log (\coth \frac{\Delta}{2}) < 1$.
Often, however, it is possible to improve this bound. For example, suppose that
$U$ is contractible in $V$ (e.g.\ if $\,V$ is simply connected) and
that $U$ is contained in a hyperbolic ball of radius $0<R<\infty$.
Lifting to the universal cover we may assume that 
$V=\DD$ and that $U=\{z\in \DD: |z|<t\}$ with
$0< t=\tanh \frac{R}{2} < 1$. The inclusion 
$(U,d_U)\hookrightarrow (\DD,d_\DD)$ has conformal derivative
$t \frac{1-|z|^2/t^2}{1-|z|^2} \leq t$ for $z\in U$. We may 
thus use $\eta=\tanh \frac{R}{2}<1$ for the contraction constant.
  Recall that for a real Birkhoff cone \cite{Bir57} one may take
  $\eta=\tanh \frac{\Delta}{4}$
 (an open interval in $\RR$ of diameter $\Delta$ is
 a ball of radius $\Delta/2$ in $\RR$).
\end{Remark}

\section{Complex Banach spaces and regularity of $\CC$-cones}
\mylabel{sec complex}
Let $X$ be a complex Banach space  and let $\calC\subset X$ be a 
$\CC$-cone (Definition \ref{def complex cone}).
We denote by  $X'$ the dual of $X$ and
we write $\la\cdot,\cdot\ra$ for the
canonical duality $X'\times X \rr \CC$.
 We will consider a bounded linear operator $T\in L(X)$ which 
 preserves $\CS$ and is
a strict and uniform contraction with respect to our gauge on $\calC$.
  We seek conditions that assure~: (1) The presence of an invariant 
  complex line (existence of an eigenvector of non-zero eigenvalue) and 
 (2) A spectral gap. In short, an invariant line appears when the 
 cone is not too `wide'
 and the spectral gap when, in addition, the cone is not too `thin'.

\begin{Definition}
\mylabel{Def aper one}
 Let $\calC\subset X$ be a $\CC$-cone in a complex Banach space
(in  section \ref{section real case}
we will use the very same
definition for
a  real cone in a real Banach space).
 When $m\in X'$ is a
 non-zero functional,  bounded on the vector space generated by $\calC$,
 we define the aperture of $\calC$ relative to $m$~:
  \[ K(\calC;m) = 
          \sup_{u\in \calC^*} 
              \frac{\|m\|\; \|u\|}{|\la m,u\ra|} \in [1,+\infty].\]
We define the aperture of $\calC$ to be~:
$\displaystyle  K(\calC) = \inf_{m\in X'^{*}} K(\calC;m) \in [1,+\infty]$.
\end{Definition}

 %Regarding the generality of the notion 
 %of `bounded sectional aperture' introduced below,
 %cf.\ Remark \ref{Remark norm-directed}.

\begin{Definition}
\mylabel{def regular cone}

\mbox{}
\begin{enumerate}
\item[(1)]
  We call $\calC$  {\em inner regular} 
 if it has non-empty  interior in $X$.

 We say that $\calC$ is {\em $T^{n_0}$-inner regular} (with $n_0\geq 0$)
  if there are $r>0$ and $x_0\in \calC$ 
 so that $x_0+T^{n_0} B(0,r) \subset \calC^*$
(when $n_0=0$ the cone is inner regular).
\item[(2)]
We say that $\calC$ is {\em outer regular} if $K(\calC)<+\infty$.

 We say that $\calC$ has {\em $K$-bounded sectional aperture}
(with $1\leq K < +\infty$) iff 
for every pair $x,y\in X$, the  sub-cone 
	 $\Span\{x,y\}\cap\; \calC$ is of
$K$-bounded aperture, i.e.\
there is a non-zero  linear functional,
$m=m_{\{x,y\}}\in \Span\{x,y\}'$,
 such that
   \begin{equation}  |\la m,u\ra| \geq  \frac1K \|u\|\, \|m\|, 
       \ \ \ \forall u\in  \Span \{x,y\} \cap \calC.
       \mylabel{m functional}
       \end{equation}
\item[(3)]
  We say that $\calC$ is {\em regular} iff the cone
   is inner and outer regular.\\
%\item[(3)]
% We call $\calC$   {\em outer regular} (in $X$)
%iff the dual cone,
% $\calC'$,  has non-empty interior in $X'$.
%\item[(4)]
%  We say that $\calC$ is {\em regular} iff the cone
%   is inner and outer regular.\\
\end{enumerate}
\end{Definition}

\begin{Remarks}
  When a cone is of bounded aperture then the cone has a bounded 
  global transverse section not containing the origin.
  This is often a too strong 
  requirement. For example, in $L^1$-spaces this is usually ok
  but not  in $L^p$ with
  $1<p\leq +\infty$ unless we are in finite dimensions.
  Being inner regular means containing an open
  ball and this typically fails in $L^p$ for $1\leq p < + \infty$,
  again with the exemption of the finite dimensional case.
  The notions of bounded sectional aperture and 
  $T^{n_0}$-inner regularity, respectively,
  are more flexible and may circumvent the
  two above-mentioned restrictions. We 
  illustrate this in Example \ref{ex real Jentzsch} and
  Theorem \ref{Complex Jentzsch}.

\end{Remarks}

It is necessary to create a passage between
the cone-gauge and the Banach space
norm. The regularity properties defined above will enable us to do so
through the following two Lemmas~:

\begin{Lemma}
\mylabel{lemma inner regular} Let $\calC$ be a $\CC$-cone and
let $x\in\CS$, $u\in X$. Suppose that there is $r>0$ such  that
$x+t u\in \CS$ 
for all $t\in\CC$ for which  $|t|<r$.
Then

\begin{enumerate}
\item[(1)]
$\ds  d_\calC (x,x+t u) \leq \frac{2}{r} |t| + o(|t|).$ 
\item[(2)] If $m$ is a linear functional on $\Span\{x,u\}$,
which never vanishes
 on the  punctured subcone, $\Span\{x,u\}\cap \calC^*$, then also~:\
$ r\; |\la m,u\ra| \leq |\la m,x\ra| $.
\item[(3)] If $\calC$ is of $K$-bounded sectional aperture then~:\
$\ds  \| u \| \leq \frac{K}{r} \|x\|$ . 
\end{enumerate}
\end{Lemma}
Proof: Let $|t|<r$. Using
	(\ref{slice}) and the scale-invariance of the cone we see that\\
$D(x,x+tu)= \{ \lambda \in \hatCC: x + \frac{1-\lambda}{2} \, t u \in \calC\}$.
Our hypothesis then implies that
$D(x,x+t u)$ contains a disc of radius
 $\ds \frac{2r}{|t|}$,
centered at 1. Shrinking a domain increases hyperbolic distances, whence
\[ d_\calC(x,x+t u) \leq d_{B(1,\frac{2r}{|t|})}(-1,1) =
         d_{\DD}(0,\frac{|t|}{r}) =
    \log \frac 
          {r + |t|}
          {r - |t|}
           =
            \frac{2}{r} |t| + o(|t|).   \]
If $m$ is non-zero on the punctured subcone,
then $0 < |\la m,x+t u\ra | = |\la m,x\ra+t \la m,u\ra| $
for all $|t|<r$ and this implies the second claim.
For the last assertion let $m$ be as in (\ref{m functional}) with $\|m\|=K$.
Possibly after multiplying $x$ and $u$
with complex phases we may assume that
$\la m,x\ra \geq  \la m, ru\ra > 0$.
Then $2r \|u\| \leq
\|x+ru\|+\|x-ru\| \leq \la m,x+ru + x-ru\ra \leq 2K \|x\|$. \Halmos
\mbox{}

\begin{Lemma}
\mylabel{bounded sectional aperture}
Let  $\calC$ be a closed complex cone of $K$-bounded sectional aperture. 
Then $\calC$ is proper, whence a $\CC$-cone (Definition \ref{def complex cone}).
 If $x,y\in \calC^*$ and  
 $m=m_{\{x,y\}}$ is a functional associated to the subcone
$\Span\{x,y\}\cap \calC$
  as in 
(\ref{m functional}) then~:
    \[ \|\frac{x}{\la m,x\ra}- \frac{y}{\la m,y\ra}\|\, \leq\, 
      \frac{4 K}{\|m\|} \tanh \frac{d_{\calC} (x,y)}{4} \, \leq \,
              K \ \frac{d_{\calC} (x,y)}{\|m\|} .\]
\end{Lemma}

Proof: We normalize the functional so that $\|m\|=K$.
Then $\|u\| \leq |\la m,u\ra| \leq K \|u\|$
 for all $u\in\Span\{x,y\}\cap \calC$.
Denote $\hatx= \frac{x}{\la m,x\ra}$
and $\haty= \frac{y}{\la m,y\ra}$ and consider,
as a function of $\lambda\in \CC$,
the point $u_\lambda=(1+\lambda) \hatx+ (1-\lambda)\haty$.
When $u_\lambda\in\calC$ the properties of $m$ show that
$\|u_\lambda\|\leq |\la m,u_\lambda\ra|=|(1+\lambda)+(1-\lambda)|\equiv 2$
and therefore,
\[
|\lambda|\; \|\hatx-\haty\|  \leq \|u_\lambda\| +
    (\|\hatx\|+\|\haty\|) \leq 4.\]
Setting $R=\frac{4}{ \|\hatx-\haty\|}\in [2,+\infty]$ we see that 
$D(\hatx,\haty)\subset \overline{B(0,R)}$. The radius $R$ is bounded iff
$x$ and $y$ are independent so the cone is proper.
Enlarging  a domain decreases hyperbolic distances so
\[ d_\calC(x,y)
       = d_{D^o(\hatx,\haty)}(-1,1) 
       \geq d_{B(0,R)}(-1,1) =
   d_\DD(
        \frac{1}{R},
        -\frac{1}{R} 
         )= 2 \log 
        \frac{1+\frac{1}{R}}
       {1 -\frac{1}{R} }.
    \] 
Therefore, $\ds \frac{\|\hatx-\haty\|}{4}=\frac{1}{R} \leq 
              \tanh \,\frac{d_{\calC} (x,y)}{4}\leq 
              \frac{d_{\calC} (x,y)}{4} $, and the  stated bound follows.
\Halmos\\
\mbox{}

\begin{Theorem}
\mylabel{invariant line}
Let $\calC$ be a $\CC$-cone of
$K$-bounded sectional aperture. Let $T\in L(X)$
be  a strict cone-contraction, i.e.\
$T: \calC^*\rightarrow\calC^*$ with
$\Delta=\diam_\calC T (\calC^*) < \infty$.
Then $\calC$ contains a unique  $T$-invariant complex
line, $\CC h$.
\end{Theorem}

Proof: Let $x_0\in\CS$ and set
$e_1=Tx_0/\|Tx_0\| \in T(\CS) \subset \CS$. We will construct
a Cauchy-sequence $(e_n)_{n\in\NNs}$ recursively.
Given $e_n$,  $n\geq 1$
choose, as in Definition \ref{def regular cone} (2),
 a functional $m_n\in X'$ normalized so that  $\|m_n\|=K$,
  associated to the subcone
  $\Span \{e_n,T e_n\}\cap \calC$.
Set $\lambda_n=\la m_n,Te_n\ra/\la m_n,e_n\ra$
(for which we have the bound $0< |\lambda_n| \leq \|T\|\;K $) and 
define the next element in our recursion~:
 \[e_{n+1} = 
      \frac{\lambda_n^{-1} T e_n }
      {\| \lambda_n^{-1} T e_n  \|} \ \in\  T^{n+1}\CS .\]
Using Lemma \ref{bounded sectional aperture} and then Lemma
\ref{lemma contraction} (with a contraction constant $\eta<1$) we obtain
for $n\geq 1$~:
   \[  \|
       \frac {e_n} {\la m_n,e_n\ra} -
       \frac {Te_n} {\la m_n,Te_n\ra} 
       \| \;\leq
       \;\dCC (e_n,Te_n)
       \;\leq\; \diam T^{n} \calC^* \leq \Delta \eta^{n-1} .\]
As $1\leq |\la m_n,e_n\ra|\leq K$ and $|\la m_n,Te_n\ra|\leq \|T\| \, K$
we get~:
\begin{equation}
   \| e_n - \lambda_n^{-1} Te_n \|
         \; \leq \; K \Delta \eta^{n-1}
   \ \ \ \mbox{and} \ \ \ \
   \| \lambda_n e_n -  Te_n \|
         \; \leq \; \|T\|\; K \; \Delta \eta^{n-1}.
   \mylabel{en bounds}
\end{equation}
 Noting that $\|e_n\|=1$,
the first inequality  
    implies~:
\begin{equation}
 \| e_n - e_{n+1} \| \leq 2 K \Delta \eta^{n-1} .
\mylabel{e bounds}
\end{equation}
The sequence, $(e_n)_{n\in\NN}$, is therefore Cauchy, whence  has a limit,
  $ h = \lim_n e_n \in \CS$, $\|h\|=1$.
The limit belongs to $\calC$ because the cone was assumed closed.
 Writing 
 $(\lambda_{n+1}-\lambda_n) e_{n+1} =
   (T-\lambda_n)e_n + (\lambda_{n+1}-T)e_{n+1} + (T-\lambda_n) (e_{n+1}-e_n)$
   and using
 the second inequality in
   (\ref{en bounds}) as well as (\ref{e bounds}) and $|\lambda_n|\leq \|T\|\;K$ we obtain
 \begin{equation}
    |\lambda_n - \lambda_{n+1}| \;\leq 
    \;(1+\eta+ (2+2K))\;\|T\|\;
       K \; \Delta\eta^{n-1}  ,
\mylabel{l bounds} \end{equation}
so  also the limit $\lambda=\lim_n \lambda_n$ exists.
But $\|Th- \lambda h\|= \lim_n 
 \|Te_n - \lambda_n e_n\|=0$  shows that $Th=\lambda h \in \CS$
which implies that  $\lambda \neq 0$, whence  that
$\CC h \subset \calC$ is a
  $T$-invariant complex line.
Suppose that also $\CC k \subset \calC$ (with $k\neq 0$) is 
  $T$-invariant.
  Then $d_\CC(h,k) \leq 
  \eta\, d_\CC(Th,Tk) =\eta\, d_\CC(h,k) \leq \eta \Delta < + \infty$ and
  this implies 
  $d_\CC(h,k)=0$ so the two vectors must be linearly dependent. 
Thus,  $\CC h$ is
unique.\Halmos\\

%\begin{Remark}
%\mylabel{rmk lower bound}
%Without further knowledge on the cone contraction there is no
%a priori lower bounds on $|\lambda^*|/ \|T\|>0$.
%The expression,
%$\inf_{u\in \CS} \|Tu\| / \|u\|$, does give a 
%lower bound for the absolute value of the eigenvalue,  but 
%this lower bound could turn out to be  zero 
%(in which case the bound is useless).
%This is not a complex phenomenon but occurs already in a real setup~:
%The reader may e.g.\ 
%consider  the projective contraction of $\RR_+^2$
% obtained through multiplication by the matrix
% $\left(  \begin{array}{cc}
%          k & 1\\
%          k^2 & k 
%          \end{array} \right)$ as  $k>0$ tends to zero.
%\end{Remark}

\begin{Theorem}
\mylabel{spectral gap}
Let $T\in L(X)$ and
let $\calC$ be a $\CC$-cone of
$K$-bounded sectional aperture which is \ {\em $T^{n_0}$-inner regular} for
some $n_0\geq 0$.
 Suppose that $T$ is a strict cone-contraction, i.e.
$T: \calC^*\rightarrow\calC^*$ with
$\Delta=\diam_\calC T (\calC^*) < \infty$.
Then $T$ has a spectral gap.
\end{Theorem}

Proof: By the previous Theorem $T$ has a unique eigenvector in the cone,
 $h\in \calC^*$,
with a non-zero complex eigenvalue, $\lambda$. 
In order to simplify the notation we replace
$T$ by  $(\lambda)^{-1}T$ and  
assume thus that there is $h\in \calC$, $\|h\|=1$
 for which $Th=h$.
 A slight complication is that $h$ need not be in the interior of $\calC$,
 or even worse, the interior of $\calC$ may be empty.
 $T^{n_0}$-inner regularity
 (Definition \ref{def regular cone} (1))
  allows us to proceed as follows:
  Let $x_0\in\calC^*$, $\|x_0\|=1$, $n_0\geq 0$ and $r>0$ be such
 that $x_0+T^{n_0} B(0,r) \subset \calC^*$.
We write  $x_n=T^{n} x_0$, $n\geq 0$ for the iterates of $x_0$.
By taking limits  in equations 
(\ref{e bounds}) and  (\ref{l bounds})
we see that
the sequences, $(e_n)_{n\in\NN} \subset \calC$ 
and $(\lambda_n)_{n\in\NN}\subset \CC^*$,  constructed in the preceding
theorem verify~:
\[
 \|e_n-h\| \leq \frac{2K\Delta}{1-\eta} \eta^{n-1} 
      \ \ \ \mbox{and} \ \ \
 |\lambda_n-1| \leq \frac{3+\eta+2K}{1-\eta}\|T\|\, 
     K\,\Delta \eta^{n-1}, \ n\geq 1 .\]

Then $\left|\|T e_n\| - 1 \right| \leq
 \| T e_{n} - Th\| \leq
  \|T\|\frac{2K\Delta}{1-\eta} \eta^{n-1}$,
so for all $n\geq 1$~:
\begin{equation}
\|x_n\| = \|x_{n-1}\|\; \|Te_{n-1}\|
      = \|Tx_0\| \prod_{k=1}^{n-1} \|Te_k\| \leq
    \|T\| \prod_{k=0}^\infty (1+ \frac{\|T\|\;2K\Delta}{1-\eta}
      \eta^k) \equiv M < \infty.
      \label{x-bound}
\end{equation}
%An explicit bound is given by  $M\leq \|T^{n_0}\| \; \exp 
 %   \left( \frac{(\|T\|+\eta) 2K\Delta}{\eta(1-\eta)^2} \right)$.
We also get that
\begin{equation} 
\| x_{n+1} - x_n \| =
   \|(T  - \lambda_n) e_n + (\lambda_n-1) e_n\| {\|x_n\|} \leq
M\,\|T\|\,K\, \frac{4+2K}{1-\eta}
\Delta\;\eta^{n-1}.
\mylabel{x ineq}
\end{equation}
%shows that $x_n$ converges to a multiple of $h$
%(since $x_n\in \CC^* e_n$ and $e_n\rr h$), so that there is
% $\gamma\in \CC^*$, $|\gamma| \leq M $ for which 
%$ \|x_n - \gamma h\| \leq \const\; \eta^n $.
Now let $u\in X$.
By our choice of $x_0$ when 
 $|t|< r/\|u\|$ then
$x_0+t\,T^{n_0} u \in \calC^*$.
 By Lemma  \ref{lemma inner regular},
 \[ \dCC(x_0,x_0+t \,T^{n_0}u) \leq 
         \frac{2 |t|}{r} \;\|u\|
          + o (\|t u\|) .\]
Applying the contraction in Lemma \ref{lemma contraction}, we get
 \[ \dCC(x_n,x_n+t\, T^{n_0+n}u) \leq 
          \frac{2 |t|}{r} \|u\| \eta^n
          + o (\|t u\| \eta^n) .\]
In order to get a norm-estimate out of this
we pick a sequence,
  $m_n$, $\|m_n\|=K$
(as in equation (\ref{m functional})),
  this time associated to  
the subcones, $\Span \{x_n,T^{n_0+n} u\} \cap \calC$.
By Lemma  \ref{bounded sectional aperture},
\[
 \| 
     \frac{x_n}{\la m_n,x_n\ra} -
     \frac{x_n + t T^{n_0+n} u}{\la m_n,x_n\ra + t \la m_n,T^{n_0+n} u\ra}
     \| \leq 
          \frac{2 |t|}{r} \|u\| \eta^n
          + o (\|t u\| \eta^n) .\]
Develop the left hand side 
in $t$, multiply by  $|\la m_n,x_n\ra|$ (which is bounded by $MK)$ 
and retain the linear term  to obtain
\begin{equation}
 \| 
     x_n \frac{\la m_n,T^{n_0+n} u\ra}{\la m_n,x_n\ra} -
      T^{n_0+n} u \|
	 \leq \frac{2 M K}{r} \eta^n \|u\| 
	    , \ n\geq 0
\mylabel{xmT ineq}
\end{equation}
valid for any $u\in X$.
Let us write 
$\alpha_n=\alpha_n(u)= \frac{\la m_n,T^{n_0+n}u\ra}{\la m_n,x_n\ra}$
for the coefficient to $x_n$.
   Since $x_0+t T^{n_0}u\in \CS$
 whenever $|t| \, \|u\| \leq r$
and $T:\CS \rr \CS$ we also have
$x_n+t T^{n_0+n} u \in \CS$ for such $t$-values. 
The second half of Lemma \ref{lemma inner regular}  then shows that
 \begin{equation}
|\alpha_n| =
    \left|\frac{\la m_n,T^{n_0+n} u\ra}{\la m_n,x_n\ra}\right|
          \leq \frac1r \|u\| ,
\mylabel{alpha ineq}
\end{equation}
uniformly in $n$. Using the identity
$x_n \alpha_n - x_{n+1}\alpha_{n+1}=
 (x_n - x_{n+1})\alpha_{n} +
 T(x_n\alpha_n -T^{n_0+n}u) +
 (
 T^{n_0+n+1}u
 -
 x_{n+1}\alpha_{n+1}
 )$.
 and the three bounds (\ref{x ineq}\ -\ \ref{alpha ineq})
we obtain for $n\geq 1$~:
\[ 
   \|x_n \alpha_n - x_{n+1}\alpha_{n+1}\| \leq
    \frac{2MK}{r}
      \left( 
        \frac{2+K}{1-\eta}\|T\|\,\Delta
	+ \eta \|T\| +  \eta^2 \right) 
        \eta^{n-1} \|u\| \equiv c_2 \; \eta^{n-1} \|u\|.\]
% with $c_2 =\frac{2M\|T\|K}{r}\,\frac{(2+K)\Delta+\eta}{1-\eta}$.
Therefore, $h\,\,c^*(u) = \lim_n x_n\, \alpha_n(u)$ exists. The limit  is 
necessarily proportional to $h$
(since $x_n\in \CC^* e_n$ and $e_n\rr h$)
 and because of (\ref{alpha ineq}) and (\ref{x-bound})
  we also have $|c^*(u)| \leq \frac{M}{r}\|u\|$.
Then,
\[ \|x_n \alpha_n - h \; c^*(u) \| \leq \frac{c_2}{1-\eta} \eta^{n-1} \|u\|,\]
so that
\[ \|h \; c^*(u) - T^{n_0+n} u\| \leq
          \left( \frac{c_2}{1-\eta} + \frac{2 M K}{r}\eta \right)
      \eta^{n-1} \|u\| \equiv C \eta^{n-1} \|u\|. \]
Linearity of $T$ implies that the mapping 
$u \rr c^*(u)=\la c^*,u\ra\in \CC $ must be  linear,
 and as a linear functional
 it is bounded in norm by $M/r$.
 % For any given $x\in \CS$,
 %$T^n x$ does not tend to zero. Thus
 %$\la c^*,x\ra\neq 0$ which shows that
 %$c^*\in \calC^{'*}$ is  an element of the dual cone.
  Finally,
this time returning to the unnormalized operator, we have shown that
 \begin{equation}
       \left\| h \la c^*,u\ra - 
       \left(\lambda^{-1}{T}\right)^{n_0+n} u \right\|
              \leq \; C \eta^{n-1}  \|u\|,
        \ \ \forall n\in\NN,\ \forall u\in X,
     \mylabel{eq projection}
  \end{equation}
with $C<+\infty$.
It follows that $\lambda$
 is a simple eigenvalue of $T$ corresponding to the
eigenprojection, $u \rr  h  \la c^*, u\ra$ and that the remainder
 has spectral radius not
exceeding $\eta |\lambda|$. \Halmos\\

%\begin{Remark}
%Some explicit estimates for the constants  in (\ref{eq projection})~:
%  \[ \|c^*\| \leq M/r \ \ \ \mbox{and} \ \ \
%     C \leq 
%       \frac{6 M K \|T\|}{r |\lambda^*|}
%         \frac{1+\Delta}{\eta(1-\eta)^2}, \]
%with $M= \frac{\|T^{n_0}\|}{|\lambda^*| }
%\exp 
%    \left( \frac{\|T\| 4K\Delta}
%     {|\lambda^*| \eta(1-\eta)^2} \right)$
% [see also Remark \ref{rmk lower bound}].
%\end{Remark}

\begin{Example}
Let $X$ be a complex Banach space and consider
$e\in X$, $\ell\in X'$ with $\la \ell,e\ra = 1$.
We write $P=e \otimes \ell$ for the associated one dimensional projection.
For
 $0<\sigma<+\infty$ we set
\begin{equation}
    \calC_\sigma= \{x \in X : \|(1-P) x\| \leq \sigma \|Px\| \}.
    \mylabel{C standard}
\end{equation}
Then $\ds B(e,\frac{\sigma \|e\|}{1+(1+\sigma)\|P\|})\subset \calC_\sigma$
and $\ds K(\calC_\sigma) \leq (1+\sigma)\|P\|$ so that
$\calC_\sigma$ is a regular $\CC$-cone.
%cf.\ Figure \ref{fig cone ex}.
 Furthermore, if $0<\sigma_1<\sigma<+\infty$  a calculation shows that
$\diam_{\calC_\sigma} \calC_{\sigma_1}^* < + \infty$.
\end{Example}

\begin{Remark}
 We have the following characterization of the 
spectral gap property~:
A bounded linear operator,
$T\in L(X)$, has a spectral gap iff it is a strict contraction of
a  regular $\CC$-cone.
Proof: One direction is the content of Theorem \ref{spectral gap}
(since a regular cone in particular is of
uniformly  bounded sectional aperture). For the other direction one uses
the spectral gap projection $P$ 
 to construct an
 adapted norm (equivalent to $\|\cdot\|$)~:
 $\|x\|_\theta = \|Px\|+ \sum_{k\geq 0} \theta^{-k} \|T^k(1-P)x \|$
 for some fixed choice of   $\theta\in (\eta,1)$. Using this norm
 to define the cone family in (\ref{C standard}),
 it is not difficult to see that
 $T$ is a strict and uniform contraction of $\calC_\sigma$, $\sigma>0$.
\end{Remark}

\section{Real cones}
\mylabel{section real case}

Let $X_\RR$ denote a real Banach space.
Recall that
a subset $\CCR\subset X_\RR$ is called a (real) proper closed convex cone
 if it is closed and convex and if

\begin{eqnarray}
 \RR_+ \CCR &=& \CCR ,\\
 \CCR \cap -\CCR &=& \{0\} .
\end{eqnarray}
We note that convexity is a useful property that {\em a fortiori} is
lost when dealing with complex cones.
In the following, we will refer to a real
 proper closed convex cone as
an {\em $\RR$-cone}. We assume throughout that such a cone is
non-trivial, i.e.\ not reduced to 
a point.
Given an $\RR$-cone one associates a projective (Hilbert) metric
for which we here give two equivalent definitions 
(for details we refer to \cite{Bir57,Bir67}).
 The first, 
 originally given by Hilbert,
  uses cross-ratios 
 and is very similar to our complex cone gauge~:
Let $\hatRR=\RR\cup\{\infty\}$ denote the extended real line
 (topologically a circle).
 For $x,y\in\CCR^*\equiv \CCR-\{0\}$,
we write 
\begin{equation}
 \ell(x,y)=\{t\in\hatRR: (1+t)x + (1-t)y\in \CCR\cup -\CCR\}
 \mylabel{def ell}
\end{equation}
with the convention that $\infty\in\ell(x,y)$ iff
$x-y\in \CCR\cup -\CCR$.
Properness of the cone implies that $\ell(x,y)=\hatRR$  iff $x$ and
$y$ are co-linear. In that case we set their distance to zero.
 Otherwise, $\ell(x,y)$ is a closed (generalized)
segment $[a,b]\subset \hatRR$ 
containing the segment $[-1;1]$,
see Figure \ref{fig cone embed} in section
 \ref{sec canonical complexification}.

The logarithm of the cross-ratio of $a,-1,1,b\in\hatRR$,
 \begin{equation}
 d_{\CCR}(x,y)  =  R(a,-1,1,b) =
          \log \frac{a-1}{a+1}  \;
           \frac{b+1}{b-1} , 
 \mylabel{Hilbert metric}
\end{equation}
then yields the Hilbert projective distance 
between $x$ and $y$. 
Birkhoff \cite{Bir57}
found an equivalent definition of this distance~: 
 For $x,y\in\CCR^*\equiv \CCR-\{0\}$, one defines
 \begin{equation}
   \beta(x,y) = \inf\{ \lambda >0: \lambda x - y \in \CCR\} \in (0,+\infty]
   \mylabel{Birkhoff beta}
 \end{equation}
 in terms of which~:
 \begin{equation}
      d_{\CCR}(x,y)= \log \left(  \beta(x,y) \beta(y,x) \right) \in [0,+\infty].
     \mylabel{Birkhoff def}
 \end{equation}
 A simple geometric argument  shows that indeed 
 the two definitions are  equivalent.
 
% \begin{Remark}
%The definition of Hilbert is really a 
%projective distance on the bi-cone
% ($\CCR \cup - \CCR$) which is $\RR$-invariant but not
% convex.
% By contrast Birkhoff's approach uses
% the convexity of the cone $\CCR$, which is only $\RR_+$-invariant. 
% Depending on the circumstances one definition
% may be easier to use than the other. 
%\end{Remark}

%As an example (which will be used below) one may define the following
% family of (real) cones in $\CC$:
%For $0\leq \theta \leq \frac{\pi}{2}$  set
%  \begin{equation}
%      \CC_\theta = \{ r e^{i\phi} : r\geq 0,
%      -\theta \leq \phi \leq \theta \}
%       \mylabel{C theta}
%   \end{equation}
%A straight-forward calculation shows that~:
%  \begin{equation}
%     \diam_{\CC_{\pi/2}}  \CC_\theta = \tan \theta/2.
%     \mylabel{C diameters}
%  \end{equation}

Given a linear functional, $m\in \XP_\RR$, the image of the cone,
$\la m,\CCR\ra$, equals either $\{0\}$, $\RR_+$, $\RR_-$ or $\RR$.
One defines the dual cone as $\CCPR=\{m\in \XP_\RR: m_{|\CCR} \geq 0\}$
and using  Mazur's Theorem, cf.\ e.g.\ \cite[p.\ 88]{Lang93}, 
 one sees that the $\RR$-cone itself may be recovered  from~:
  \begin{equation}
     \CCR = \{ x\in X_\RR :  \la m,x\ra \geq 0 , \ \forall m\in\CCPR \}.
     \mylabel{dual cone char}
     \end{equation}
\mbox{}

Given an $\RR$-cone $\CCR$ we  use 
   Definition \ref{Def aper one} (replacing $\CC$ by $\RR$, complex
   by real) to define
  the aperture of $\CCR$. It is given as
   the infimum of $K$-values for which there exists
a linear functional $m$ satisfying (see Figure \ref{real cone figure})
  \begin{equation}
   \|u\| \leq \la m,u \ra \leq K\|u\|, \ \ \  u\in \CCR.
    \mylabel{R aper}
  \end{equation}
     
\begin{figure}[htbp]
\begin{center}
{\input{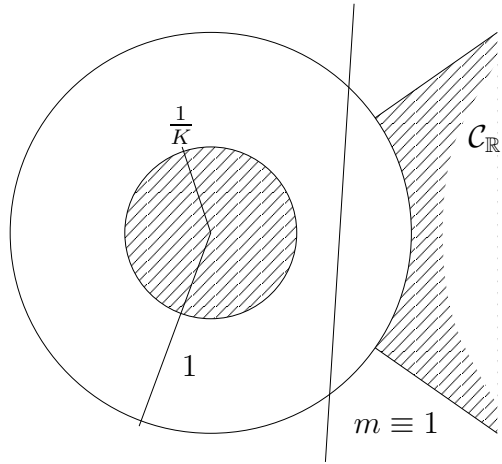}}
\caption{A real cone $\CCR$ of $K$-bounded sectional aperture}
\mylabel{real cone figure}
\end{center}
\end{figure}

\begin{Lemma}
\mylabel{L real aper}
The aperture, $K(\CCR)\in [1,+\infty]$,   of an $\RR$-cone,
 $\CCR\subset X_\RR$, is determined by
  \begin{equation}
    \frac{1}{K(\CCR)} 
       \;=\; \inf\;  \{\; \frac{\|x_1+ \cdots+ x_n\| } 
              {\| x_1 \| + \cdots+ \| x_n\|}\; :\; \
	         x_i\in \CCSR, n\geq 1\} .
      \mylabel{real aper}
\end{equation}
\end{Lemma}
Proof: 
Let $x_1,\ldots,x_n\in\CCSR$  and note that
$a=\sum_1^n x_i / \sum_1^n \|x_i\|$ belongs to
 $A \equiv \Conv (\CCR \cap \partial B(0,1))$, the convex hull
of cone-elements of norm one.
The reciprocal of the right hand side in (\ref{real aper}) 
therefore equals 
    $r = \inf \{ {\|a\|} : a \in A\}\in [0,1]$.
 Suppose that $r>0$.
 Then
${B(0,r)}$ and $\Cl {A}$ are disjoint convex subsets. The vector
 difference,
$Z=\{a-b: a\in\Cl A , b\in  B(0,r)\}$, is open, convex and does not contain the origin,
whence  \cite[Lemma 2.2, p.89]{Lang93} there is  $\ell\in X'_\RR$ 
whose kernel does not intersect $Z$. We may normalize $\ell$ so that
 \[ {B(0,r)} \subset \{\ell<1\}  \ \ \
      \mbox{and} \ \ \ \Cl {A} \subset \{\ell \geq 1\}. \]
Then
  \begin{equation}
   \| x\| \leq \la \ell,x\ra \leq  \frac{\|x\|}{r} , \ \ \ \forall x\in \CCR^*,
  \mylabel{l-inequality}
  \end{equation}
and therefore $K(\CCR)\leq \frac1r$. 
To get the converse inequality let  $m$ be
positive and verify (\ref{R aper}). Then
$\sum \|x_i\| \leq 
\sum \la m,x_i \ra =
 \la m,\sum x_i \ra \leq
K \| \sum x_i\|.$ \Halmos\\

\begin{Lemma}
\mylabel{CR bd aperture}
Let $\CCR\subset X_\RR$ be a $d$-dimensional
 $\RR$-cone of $K$-bounded sectional aperture.
Then $\CCR$ itself is of
$dK$-bounded
aperture.
\end{Lemma}
Proof:
Let $F\subset \RR^d$. By a theorem of Caratheodory, 
a point in the convex hull, $\Conv F$,
  is {\em a fortiori}
in the convex hull of $d+1$ points in $F$
(see e.g.\ \cite[p.73]{Rud91}). If 
$x\in \partial\; \Conv F$,  we may even write it as a limit of
convex combinations of $d$ points in $F$. Now, apply this to
the set $A$ in the proof of the previous Lemma.
In the formula, (\ref{real aper}) it thus suffices to
consider $d$ cone-elements 
 which we may  order decreasingly
according to their norm,
 $\|x_1\| \geq \|x_2\| \geq \cdots \geq \|x_d\|$. 
Using Lemma 
\ref{L real aper} with $n=2$,
 the K-bounded
 sectional aperture implies that
     \[
     \|x_1+\cdots+x_d\| \geq
            \frac{1}{K} \|x_1\| + \frac{1}{K} \|x_2+\cdots+x_d\| \geq 
	    \frac{1}{K} \|x_1\| \geq 
  \frac{1}{K} 
\frac{\|x_1\|+\cdots+ \|x_d\|}{d}.
\]
Thus
\[
\frac{\|x_1+\cdots+x_d\|}
{\|x_1\|+\cdots+ \|x_d\|}
\geq \frac{1}{dK}  ,
\]
 %\frac{1}{d} \left(\frac{1}{K} + \frac{1}{K^2} +
            %\cdots + \frac{2}{K^{d-1}}\right) \geq \frac{1}{dK}  .\]
 and in view of  Lemma \ref{L real aper}, we see
 that $\CCR$ is of $dK$ bounded aperture.\Halmos\\

\begin{Lemma}
\mylabel{R outer reg}
Let $\CCR\subset \RR^d$ be an $\RR$-cone. Then $\CCR$ is outer regular.
\end{Lemma}
Proof~: As in the previous Lemma it suffices to look at the supremum in
(\ref{real aper}) over $d$-tuples.
The set $A\equiv \{x_1,\ldots,x_d \in \CCR : \|x_1\|+\cdots+\|x_d\|=1\}$
is compact and $\|x_1+\cdots+x_d\|$ is continuous and non-vanishing on $A$,
whence has a minimum, $r>0$.
It follows that $K(\CCR)\leq \frac{1}{r}<+\infty$.\Halmos\\

\begin{Remark}
\mylabel{Remark norm-directed}
In the literature an $\RR$-cone $\CCR$ is 
said to be norm-directed (with a constant $0<K<\infty$) if 
$\|x-y\|\leq K \|x+y\|$, $\forall x,y\in\CCR$.
For an $\RR$-cone
our notion of uniformly bounded sectional aperture is equivalent
(up to a small unavoidable loss in constants) to that of being  norm-directed.
 To see this note that if $\CCR$ is of $K$-bounded sectional aperture 
 and $\ell$ verifies (\ref{l-inequality}) then
    $\forall x,y\in \CCR$~:
 \[ \|x-y\| \leq \|x\|+\|y\| \leq
  \langle \ell, x\rangle + \la \ell,y\ra
   = 
  \langle \ell, x+y\rangle \leq K \|x+y\|, \]
 which shows that $\CCR$ is $K$-norm-directed. 
 Conversely, if $\CCR$ is $K$-norm-directed then
\[ \|x\| + \|y\| \leq \|x+y\| + \|x-y\| \leq
    (1+K) \|x+y\| \]
and  Lemma \ref{L real aper} shows that $\CCR$ is of $(K+1)$-bounded
sectional aperture.
For example, $(\RR^d_+,\|\cdot\|_1)$ is 1-norm directed and of
1-bounded aperture, whereas $(\RR^d_+,\|\cdot\|_\infty)$ is 
1-norm directed, of 2-bounded sectional aperture
 but only of $d$ bounded  aperture.
 Lemma \ref{lemma inner regular} (3)
 is a complex cone-analogue of being norm-directed.
 \end{Remark}

\begin{Theorem} \mylabel{Thm Birkhoff}
Let $A\in L(X_\RR)$ and let $\CCR\subset X_\RR$ be an $\RR$-cone
which is $A^{n_0}$-inner regular
(for some $n_0\geq 0$) and $K$-norm-directed. Suppose that
$A$ is a strict cone-contraction, i.e.\
$A:\CCSR \rr \CCSR$ with $\Delta_A=\diam_{\CCR} A(\CCR^*) < +\infty$.
Then $A$ has a spectral gap.
More precisely, there is $\lambda>0$ and a one dimensional
projection $P$ for which $\lambda^{-1} A-P$ has spectral radius
not greater  than $\tanh\frac{\Delta_A}{4} < 1$.

Proof: See 
\cite{Bir57,Bir67} in the case of
$\CCR$ being inner regular and $K$-norm directed.
When $\CCR$ is assumed only $A^{n_0}$-inner regular 
with $n_0>0$ one may 
either adapt the proof of Birkhoff (easy)
or use Remark \ref{C contains R} below.
\end{Theorem}

\begin{Corollary}
\mylabel{R finite dim}
Let $\CCR$ be an $\RR$-cone in $\RR^d$, $d<+\infty$
 and suppose that 
$A\in L(X_\RR)$ verifies $A(\CCSR) \subset \Int \CCR$. Then 
$A$ has a spectral gap.
\end{Corollary}
Proof: Implicitly it is assumed that $\CCR$ has non-empty interior.
Lemma \ref{R outer reg} shows that $\CCR$ is outer regular,
in particular, norm-directed.
Local compactness of $\RR^d$
 implies that $\diam_{\CCR} A(\CCR^*) <+\infty$ so the Corollary
follows from Birkhoff's Theorem.\Halmos\\

Let $\CCR\subset X_\RR$ be an $\RR$-cone. It is standard to 
write $x\preceq y  \Leftrightarrow y-x\in \CCR$ for the induced partial
ordering of
$x,y\in X_\RR$.
For $A,B\in L(X_\RR)$, we also write 
$A\preceq B  \Leftrightarrow \forall x\in\CCR \; : \;
A(x)\preceq B(x)$. 

The following dominated cone contraction theorem
is trivial in the context of an $\RR$-cone contraction.
In section \ref{C dominated} we show
 that a similar (non-trivial)
result holds in the complex case.

\begin{Theorem}
Let $A,P\in L(X_\RR)$ be contractions of the $\RR$-cone $\CCR$.
Suppose that there are constants $0< \alpha \leq \beta < +\infty$
for which
      $\alpha P \preceq A \preceq \beta P$.
Then
   \[ \diam_{\CCR} A(\CCR^*) \leq 2 \log \frac{\beta}{\alpha}
           + \diam_\CCR P(\CCR^*).\]
	\mylabel{real domin gap}   
\end{Theorem}
Proof: Given $x,y\in \CCR^*$, suppose that
$\lambda,\lambda'>0$ are such that
$\lambda Px-Py \in \CCR$, $\lambda'Py-Px\in\CCR$.
Then also $\lambda\beta Ax -\alpha Ay \in\CCR$
and $\lambda'\beta Ay -\alpha A x\in\CCR$ so that
\[ d_{\CCR} (Ax,Ay) \leq 2\log \frac{\beta}{\alpha} +
 \log(\lambda\lambda')\] and the claim follows
 by Birkhoff's characterization
 (\ref{Birkhoff def}).\Halmos\\
 
 \begin{Example}
    \mylabel{ex real PF}
    \mylabel{ex real Jentzsch}
 \mbox{}\\[-5mm]
 \begin{enumerate}
   \item [(1)]
   Let $A\in M_n(\RR)$ and suppose that
 $0< \alpha\leq A_{ij} \leq \beta <+\infty$ for all indices.
 Setting $P_{ij}\equiv 1$ we see that
  \[ P((\RR_+^n)^*) = \{(t,\ldots,t): t>0\}\]
  which is of zero projective diameter in $\RR_+^n$.
  By Theorem \ref{real domin gap}
  we recover the standard result~:
    \[ \diam_{\RR^n_+} A((\RR_+^n)^*) \leq
     \Delta_A= 2 \log \frac{\beta}{\alpha}.\]
  In particular, if $\lambda_1>0$ and $|\lambda_2|$ 
  denote the leading eigenvalue and the second largest eigenvalue
  (in absolute value), respectively, then 
   $\frac{|\lambda_2|}{\lambda_1 } \leq 
   \tanh \frac{\Delta_A}{4} = \frac{\beta-\alpha}{\beta+\alpha}$.
 \item  [(2)]
 The standard Perron-Frobenius Theorem
   generalizes easily to integral operators,
   cf.\ Jentzsch's Theorem \cite{Jen12} and the generalization
   given by Birkhoff in \cite{Bir57}). We present a
   somewhat different generalization~:
   Let $(\Omega,\mu)$ be a measure space and let
   $X_\RR=L^p\equiv L^p(\Omega,\mu)$, $1\leq p \leq +\infty$.
    Let $h\in L^p_+$ ($h>0$, a.e.)
    and $m\in L^q_+$ ($m>0$, a.e.) with $q=p/(p-1)\in[1,+\infty]$ being the
    conjugated exponent so that $0<\int_\Omega h\; m\; d\mu <+\infty$.
    Let $k_A:\Omega\times \Omega\rr \RR_+$ 
    be a $\mu\otimes\mu$-measurable map.
    We suppose there are constants   $0<\alpha\leq \beta < +\infty$ 
  so that for $\mu$-almost all $x,y\in \Omega$~:
          \[ \alpha\; h(x)m(y) \leq k_A(x,y) \leq \beta \;h(x)m(y). \]
   Let $A\in L(X_\RR)$ be the integral operator defined by
   $A\phi(x)= \int_\Omega k_A(x,y) \phi(y) \, d\mu(y)$.
  Then $A$ has a spectral gap (again with a contraction rate given
   by $\frac{\beta-\alpha}{\beta+\alpha}$).

   \ \  Proof: We write 
   $\CCR=L^p_+(\Omega,\mu)$ for the cone of positive
   $L^p$-functions ($\phi\geq 0$, a.e.) 
   and compare the operator $A$ with the one-dimensional
   projection
   $ P\phi = h \; \int_\Omega m\,\phi\, d\mu$.
   Our assumption, $\int h\,m\,d\mu >0$ shows that 
   $P:\CCSR \rr \CCSR$ and 
   that  $\Delta_P=0$.
  Thus, \ $\diam_{\CCR} A(\CCR^*) \leq 2\;\log\, \frac{\beta}{\alpha}$.

   In general, $\CCR$ need not be inner regular (unless $p=\infty$).
   It is, however, $A^{n_0}$-inner regular (with $n_0=1$) 
    for any value
   of $1\leq p\leq \infty$. One
   has~: $h+Au \geq h(1-\beta\int m\,|u|\,d\mu)$ so
   $h+Au\in \CCSR$ when $\|u\|_p < \frac{1}{\beta\, \|m\|_q}$.
   To see that $\CCR$ is of 
  uniformly bounded sectional aperture let $f,g\in\CCR$ be of unit norm
  in $L^p_+$, $1\leq p<+\infty$
   and pick $\tf,\tg \in L^q_+(\Omega)$ with $\|\tf\|_q=\|\tg\|_q=1$
   (the case $p=\infty$, $q=1$ should be treated separately; we 
   leave this to the reader) and 
  $\int \tf\, f \, d\mu =
  \int \tg\, g \, d\mu = 1$. The functional
  $m(u)=\int (\tf + \tg)u \, d\mu$ then verifies
  $\|u\|_p  \leq m(u) \leq 2 \|u\|_p$ for all 
  $u\in \CCR \cap \Span\{f,g\}$. Thus $\CCR$ is of 2-bounded
  sectional aperture.
  \end{enumerate}
 \end{Example}
 
\section{The canonical complexification of a real Birkhoff cone}
\mylabel{sec canonical complexification}
A complex cone yields a genuine extension/generalization
 of the cone contraction described by Birkhoff \cite{Bir67}.
 More precisely, we will show that any Birkhoff cone may
be isometrically embedded in a complex cone, enjoying qualitatively
the same contraction properties. 

Let $X_\RR$  be a Banach space over the reals.
A complexification $X_\CC$ of $X_\RR$ is a complex Banach space,
equipped with a bounded
 anti-linear complex  involution, $J:X_\CC\rr X_\CC$,
$J^2={\bfone}$, $J(\lambda x)=\overline{\lambda} 
J(x)$, $J(x+y)=J(x)+J(y)$, $\lambda\in\CC$,
$x,y\in X_\CC$, for which $X_\RR=\half(\bfone+J) X_\CC\,$  is the real part.
 Then $X_\CC=X_\RR \oplus i X_\RR$
 is a direct sum.
[Note that this is not the same as
 regarding $X_\CC$ as a real Banach space.
For example, $\CC^n$ is a complexification of $\RR^n$ for any
$\ell^p$-norm, $1\leq p\leq \infty$, while 
the real dimension of $\CC^n$ is $2n$].
Usually $J$ will be an isometry on $X_\CC$ in which case
 the canonical projections,
$\Re=\half(\bfone+J)$ and
$\Im=\frac{1}{2i}(\bfone-J)$, have norm one.
We note that any real Banach space, $(X_\RR,\|\cdot\|_\RR)$,
admits a complexification,
 $X_\CC=X_\RR \oplus i X_\RR$ as follows~: We adopt 
 the obvious rules for multiplying by  complex numbers,
 set $J(x+iy)=x-iy$ and introduce
a norm e.g.\ using real functionals,
 \[  \|x+iy\|_{\CC} =
 \sup \{ |\la \ell,x\ra+ i \la \ell,y\ra|: \ell\in X^{'}_\RR,
         \|\ell\|_\RR \leq 1 \}.
\]
The latter norm is equivalent (within a factor of 2) to any other
norm on $X_\CC$ having as real part the given space $(X_\RR,\|\cdot\|_\RR)$.
For the rest of this section  $X_\RR$ will denote the real part of a
complex Banach space $X_\CC$.
A real linear functional, $m\in\XP_\RR$, extends to a complex linear
functional by setting $\la m,x+iy\ra=\la m,x\ra+i \la m,y\ra$ for 
$x+iy \in X_\RR \oplus i X_\RR$.
\begin{Definition}
Given an $\RR$-cone
 $\CCR\subset X_\RR$ 
we define its {\it canonical complexification}~:
  \begin{equation}
     \CCC = \{ u\in X_\CC:
           \Re \ \la m,u\ra \, \overline{\la \ell,u\ra} \, \geq \, 0,
	           \ \ \forall m,\ell \in \CCPR \}.
          \mylabel{complexification}
  \end{equation}
\mylabel{cone complexification}
\end{Definition}

%  \begin{equation}
%     \CCC = \{ x\in \XP_\CC:
%           \lambda x + \mu y : x,y\in \CCR, \lambda,\mu\in \CC,
%           |\lambda-\mu| \leq |\lambda + \mu| \}.
%  \end{equation}

\begin{Proposition}
We have the following 
{\bf polarization} identity~:
 %\begin{equation}
 %\calC_\CC =
 %       \{ \lambda v: \lambda \in \CC^*, v\in X_\CC :
 %                       \la m,v\ra \in \CC_{\pi/4},
 %      \forall m\in\CCPR 
 %\} 
 % \end{equation} 
 \begin{equation}
 \CCSC
=\{ \lambda (x+ i y) : \lambda\in\CC,
                      x\pm y \in  \CCSR \} 
    \mylabel{equiv C}
 \end{equation}
\mylabel{C cone char}
\end{Proposition}
Proof: Let $u\in\CCSC$.
Our defining condition (\ref{complexification})
 means that $\la m,u\ra$  (assume here it is non-zero) must have
 an argument that
vary within a $\pi/2$ angle as $m\in\CCPR$ varies.
 Normalizing appropriately, we may write
 $u=\lambda v$ with $\lambda\in\CC^*$ and $\Arg \la m,v\ra \leq \pi/4$.
 If we set $v=x+iy$ with $x,y\in X_\RR$ then
  $|\la m,y\ra| \leq \la m, x \ra$ for all $m\in\CCPR$.
Hence, $\la m, x\pm y\ra \geq 0$ for all such functionals and 
by (\ref{dual cone char}) this is
equivalent to $x\pm y \in\CCR$.
If $x=y$ (or $x=-y$) then we may write 
$u=\lambda(1+i) \; x$ (or
$u=\lambda(1-i) \; x$) so we may always assume $x\pm y\in \CCSR$.
\Halmos\\

\begin{Lemma}
\mylabel{C aperture}
Let $\CCR$ be an $\RR$-cone of $K$-bounded aperture. Then its 
canonical complexification,
$\CCC$, is of $2\sqrt{2}\,K$-bounded aperture.
\end{Lemma}
Proof: Let $\ell\in X'_\RR$ satisfy
   $\|x\| \leq \la \ell,x \ra \leq K\|x\|, \  x\in \CCR$ and
   extend $\ell$ to a complex linear functional. 
When $u\in \CCC$ we use polarization, Lemma \ref{C cone char}, to write
$u=\lambda (x+iy)$ with $\la \ell, x\pm y\ra \geq 0$.
Then $\|x\pm y\| \leq \la \ell, x\ra \pm \la \ell , y\ra \leq
K \|x\pm y\|$, from which $\|x\| \leq \la \ell,x\ra$ and 
$\|y\|\leq \la\ell,x\ra$
so that
$\frac12 \|x+i y\| \leq \la \ell,x\ra \leq  
  | \la \ell, x\ra +i \la \ell,y \ra|$.
  As $|\la \ell,y\ra| \leq \la \ell,x\ra$ we also have
$|\la \ell,x+iy\ra |\leq \sqrt{2}\, \la \ell,x\ra 
 \leq \sqrt{2}\, K \|x\| \leq \sqrt{2}\, K \|x+iy\|$. Therefore,
   \[\frac12 \|u\| \leq |\la \ell,u\ra| \leq \sqrt{2}\, K \; \|u\|.\Halmos \]

\begin{Proposition}
\mylabel{C regular cone}
Let $\CCR$ be an $\RR$-cone. 
If $\CCR$ is (1) inner regular/ (2) outer regular/ (3) of bounded
sectional aperture then so is its canonical complexification.
\end{Proposition}
Proof: (1) When $\CCR$ contains a ball $B_{X_\RR}(h,r)$ one verifies that
$\CCC$ contains the ball $B_{X_\CC}(h,r/2)$.\\

(2) 
As shown in Lemma \ref{C aperture},
 if $V_\RR\subset X_\RR$ is of $K$-bounded aperture then 
$V_\CC$ is of $2\sqrt{2}\,K$-bounded aperture.\\

(3) Let $u_1,u_2\in \CCSC$ and write 
$W=\Span_\CC\{u_1,u_2\} \cap \CCC$ for the sub-cone generated by
these two elements.
%$u_1=\lambda_1(x_1+iy_1)$,
%$u_2=\lambda_2(x_2+iy_2)$ as in  (\ref{equiv C}).
We also write $F=\Span_\RR \{\Re \,u_1, \Im \,u_1, \Re \,u_2, \Im \,u_2 \}$
and $V_\RR=F\cap \CCR$ which is an at most 4
and at least 1-dimensional 
$\RR$-subcone of $\CCR$. Now, if $w\in W$ then $w=\lambda'(x'+iy')$ with
$x'\pm y'\in \CCR$ and clearly also $x',y'\in F$. But then
$x'\pm y'\in V_\RR$ so that also $w\in V_\CC$,
with $V_\CC$ being the complexification of $V_\RR$.
By Lemma \ref{CR bd aperture}, $V_\CC$ is of $4\,K$ bounded
aperture so by Lemma \ref{C aperture},
$V_\CC$ and therefore also $W$
are  of $8\sqrt{2}\,K$ bounded (complex)  aperture.\Halmos\\

\begin{Theorem}
 \mylabel{thm embedding}
  Let $\CCR$ be an $\RR$-cone and let
  $\CCC$ denote its canonical complexification (\ref{complexification}).
  Then $\CCC$ is a $\CC$-cone 
 (Definition \ref{def complex cone}).
Writing $d_\CCC$ for our
projective gauge on the complex cone,
  the natural inclusion,
  \[ 
      \left( \CCSR,\dCCR \right) 
        \hookrightarrow
      \left( \CCSC,\dCCC \right) , \]
  is an isometric embedding.
\end{Theorem}

Proof:
 The set $\CCC$ is clearly $\CC$-invariant.
 Consider independent vectors, $x,y\in \CCR^*$.
 By Lemma \ref{R outer reg} any finite dimensional subcone of $\CCR$
 is outer regular, so in particular of uniformly bounded sectional aperture.
 Our previous Lemma shows that the corresponding complex cone
 is of bounded sectional (complex) aperture. But then $\CCC$ must
 be proper by Lemma \ref{lemma inner regular}(3).
 
 Regarding the embedding we may normalize  the points so that
$\ell(x,y)=[a,b]$ 
is a bounded segment in $\RR$. Define the `boundary' points,
$x_0=(1+a)x+(1-a)y$ and $y_0=(1+b)x+(1-b)y$. For any $\epsilon>0$ the point
$-\epsilon x_0 + (1+\epsilon) y_0$ is outside the closed convex
cone $\CCR$. By Mazur's Theorem, \cite[p.88]{Lang93}, we may separate
this point from $\CCR$ by a functional $\ell\in\CCPR$.
For any $\epsilon>0$ we may then find $m,\ell\in\CCPR$ 
for which
\[ \la m,x_0\ra=\la \ell,y_0\ra=\epsilon\ \ \mbox{and} \ \ 
   \la m,y_0\ra=\la \ell,x_0\ra=1 .\]
%[Examples, e.g.\ in $\RR^3$, show that one may not 
%necessarily achieve $\epsilon=0$ here].
Then $u=\mu x_0+\lambda y_0\in \CCC$ only if
$\Re \,
 (\epsilon \mu + \lambda)(\epsilon \barla + \barmu) \geq 0$
for any $0<\epsilon \leq 1$ whence only if
 \[ \Re \ \lambda\; \barmu \geq 0
  \ \ \Leftrightarrow\ \  |\lambda+\mu|^2 \geq |\lambda-\mu|^2 .\]
Conversely, when $\Re\, \lambda \, \barmu \geq 0$ and $m,\ell\in \CCPR$
then 
 \[  \Re \, \la m,u\ra \, \overline{\la \ell,u\ra} \, \geq
        \,\Re\, ( \lambda \, \barmu  )
          \,\left( 
	   \,\la m, y_0 \ra  \,\la \ell,x_0 \ra
	   \,+ \, \la m, x_0 \ra  \,\la \ell,y_0 \ra  \,\right) \, \geq \, 0 ,\]
so this condition is also sufficient. We thus have~:
$D(x_0,y_0)=\overline{\DD}$. Therefore
 $D=D(x,y)\subset \hatCC$ is a generalized disc, symmetric under complex
conjugation for which $\ell(x,y)=D(x,y)\cap \hatRR$
(see Figure \ref{fig cone embed}).
In this situation we know explicit formulae for
both  (\ref{Hilbert metric}) the real and  (\ref{C metric}) the complex
hyperbolic metrics in terms of cross-ratios so
 we get
 $\dCCC(x,y)=d_{D^o(x,y)}(-1,1)=\log R(a,-1,1,b) = \dCCR(x,y)$.\Halmos\\

\begin{figure}
\begin{center}
\epsfig{figure=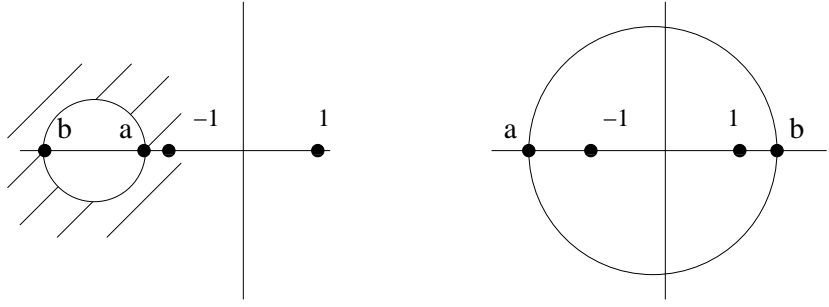,height=4cm}
\end{center}
\caption{Illustration of two possible configurations of $D(x,y)$}
\mylabel{fig cone embed}
\end{figure}

\begin{Corollary}
\mylabel{Corollary Cplus}
For $n\geq 1$ the set 
\[ \CCone = 
    \{u\in \CC^n : \Re \; u_i \overline{u_j} \geq 0, \ \forall i,j\}
    = \{u\in \CC^n : |u_i+u_j| \geq |u_i-u_j|, \ \forall i,j\}
\]
is a regular $\CC$-cone.
The inclusion, $((\RR^n_+)^*,d_{\RR^n_+}) \hookrightarrow
((\CCone)^*,d_{\CCone})$ is
an isometric embedding.
\end{Corollary}

Proof: Let $\ell_i\in (\RR^n)'$, $i=1,\ldots,n$ denote the canonical
coordinate projections then
 $\RR_+^n = \{x\in \RR^n : \la \ell_i,x\ra \geq 0,\ \forall\ i\}$ and
 $\CC_+^n = \{u\in\CC^n : \Re \;\la \ell_i,u\ra 
        \overline{\la \ell_j,u \ra}\geq 0,\ \forall\ i,j\}$.
Thus, $\CC_+^n$ is the canonical complexification 
of the standard 
real cone $\RR_+^n$.
 See Figure \ref{fig can complex} for an illustration of $\CC^2_+$.
 \Halmos\\

\begin{figure}[htbp]
\begin{center}
\input{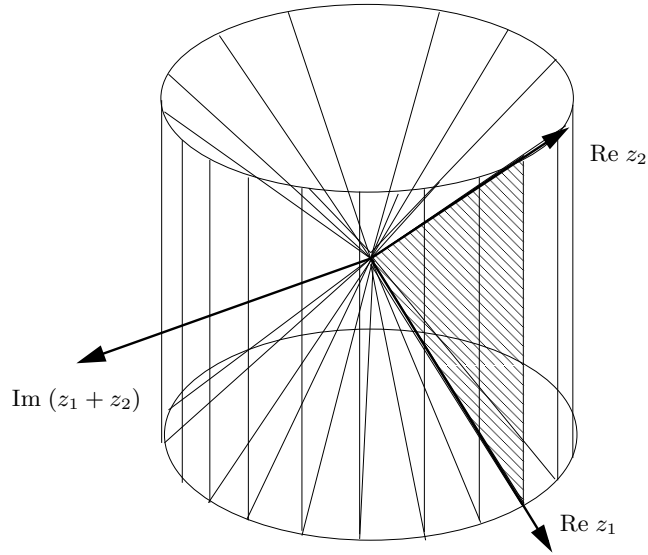}
\end{center}
\caption{An attempt to illustrate the canonical complexification $\CC^2_+$
    of $\RR_+^2$ in the coordinate system 
    $\ds (\Re\; z_1$, $\Re\; z_2$, $\Im\; (z_1+z_2))$ and
    setting  $\Im\;(z_1-z_2)=0$. We show the part of the cone
    contained in the region
     $(\Re\;(z_1+z_2))^2 + 
     (\Im\;(z_1+z_2))^2 \leq 1 $
     The shaded region shows the intersection with the real cone,
      $\RR_+^2$.
         \mylabel{fig can complex} }
\end{figure}

Below we shall need the following complex polarization
\begin{Lemma}
\mylabel{complex polarization}
Let $x\pm y \in \CCSR$ be at a distance
$\Delta=d_{\CCR}(x-y,x+y) < + \infty$.
We may then find $\alpha\in\RR$ so that
the vector $x'+iy'=e^{i\alpha} (x+iy)$, or equivalently~:
\begin{eqnarray*}
    x' &=&  x\; \cos \alpha - y \;\sin\alpha\\
    y' &=&  x \;\sin\alpha  + y\; \cos\alpha
\end{eqnarray*}
verifies:
\[  x'+ty'\in \CCR, \ \ \forall\ |t| \leq \coth \frac{\Delta}{4}.\]
In particular, $x'\pm y'\in \CCSR$ and
\begin{equation}  |\la \ell,y'\ra| \leq
     \left( \tanh \frac{\Delta}{4} \right)\; \la \ell,x'\ra\ ,
   \ \ \ \forall \ell \in \CCPR.
   \mylabel{rotated polarization}
\end{equation}
\end{Lemma}
Proof:
Possibly after replacing $(x+iy)$ by
$e^{i\phi}(x'+iy')$ for a suitable $\phi\in\RR$ we may assume that
$x\in\CCR$ and $y\notin \CCR\cup - \CCR$.
Then 
\begin{eqnarray*} \ell(x-y,x+y) 
&=& \{ t\in \hatRR: (1+t)(x-y)+(1-t)(x+y)\in \CCR\cup-\CCR\} \\
&=& \{ t\in \RR: x-ty\in \CCR \}
\end{eqnarray*}
is a real segment $[t_1,t_2]$ for which $]t_1,t_2[ \supset [-1,1]$
(see Figure \ref{fig geom}). Now write $t=\tan(\theta)$,
$\theta\in I\equiv ]-\pi/2,\pi/2[$. Since $\cos(\theta)>0$ we get
in the $\theta$-coordinate~;
\begin{eqnarray} \Theta(x,y) = [\theta_1,\theta_2]
&=& \{ \theta \in I:  x - \tan \theta \; y \in \CCR\} \nonumber \\
&=& \{ \theta \in I: \cos \theta \; x - \sin \theta \; y \in \CCR\} \nonumber \\
&=& \{ \theta \in I: \Re\; e^{i\theta} (x+iy) \in \CCR\} \label{thetaxy}
\end{eqnarray}
% If we define the function $f(\theta)=\frac{\tan(\theta)+1}{\tan(\theta)-1}$
%then the projective distance between
%$x'-y'$ and $x'+y'$ is given by
%  \[ d_(\theta_1,\theta_2) = \log \left( f(\theta_2) / f(\theta_1) \right).\]
In $\theta$-coordinates the projective distance between
$x-y$ and $x+y$ is then given by
  \[ d_(\theta_1,\theta_2) \equiv  
   \log \left( 
   \frac{\tan(\theta_2)+1}{\tan(\theta_2)-1}
      \; \times \; 
   \frac{1-\tan(\theta_1)}{1+\tan(\theta_1)}
        \right).\]

 If we do a complex rotation, $x'+iy'=e^{i\alpha}(x+iy)$ 
 with $\alpha \in J\equiv ]-\frac{\pi}2
   +\theta_1 ,\frac{\pi}{2} - \theta_2[$ then
 the last expression in (\ref{thetaxy}) shows that
 $\Theta(x',y')=[\theta_1-\alpha,\theta_2-\alpha]$.
Now, the derivative of 
$\alpha\in J \rr
  d_\CCR(x'-y',x'+y')= d(\theta_1-\alpha, \theta_2-\alpha)$
 equals 
 \[  \frac{2}{\cos(2(\theta_1-\alpha))}
   -\frac{2}{\cos(2(\theta_2-\alpha))},\]
so the minimal distance between $x'-y'$ and $x'+y'$
is obtained for
$\alpha=-(\theta_2+\theta_1)/2$. This corresponds to a symmetric
configuration in which 
$\ell(x'-y',x'+y')=[-L,L]$ with
$L=\tan\frac{\theta_2-\theta_1}{2} >1$ and
$\Delta\geq d(x'-y',x'+y') = 2 \, \log \, \frac{L+1}{L-1}$
or equivalently,
\[ L \geq \coth \frac{\Delta}{4}.\]
Since $x'+ty'\in \CCR$ whenever $|t|\leq L$ we obtain the first claim.
The second is a consequence of 
$\la \ell,x'+t y'\ra \geq 0$ for all $-L\leq t \leq L$.
\Halmos

\begin{figure}
\begin{center}
{\input{geom.1.pstex_t}}
\end{center}
\caption{ The subcone $\CCR\cap\Span\{x,y\}$ viewed in the
           $x$,$y$-coordinate system.
         \mylabel{fig geom} }
\end{figure}

\begin{Lemma}
\mylabel{real rescaling}
Let $x_1,x_2\in \CCR$ be at a distance $\Delta=d_\CCR(x_1,x_2)<+\infty$.
Through a positive real rescaling, e.g.\ replacing \ $x_1$ by $tx_1$
for a suitable  $t>0$, we may assure that
\[ |\la \ell,x_1-x_2\ra| \leq
     \left(\tanh\frac{\Delta}{4} \right)\; \la \ell,x_1+x_2\ra
    , \ \ \ \forall 
           \ell \in \CCPR.\]
\end{Lemma}
Proof: From the Birkhoff characterization (\ref{Birkhoff def}) of
 the projective distance we may rescale, say  $x_1$, to obtain
 $e^{\Delta/2} x_1-x_2\in \CCR$
 and $e^{\Delta/2} x_2-x_1\in \CCR$. Then
    $ \la \ell,x_1\ra \leq e^{\Delta/2} \la \ell,x_2\ra$ and
    $\la \ell,x_2\ra \leq e^{\Delta/2} \la \ell,x_1\ra$.
         From this we get~:
    $(e^{\Delta/2}-1)(\la \ell,x_1\ra + \la \ell,x_2\ra)-
     (e^{\Delta/2}+1)(\la \ell,x_1\ra - \la \ell,x_2\ra)=
    2 (e^{\Delta/2} \la \ell,x_2\ra
    - \la \ell,x_1\ra) \geq 0$ and similarly with $x_1$ and $x_2$
    interchanged.
    Rearranging terms the claim follows.
 \Halmos\\
 
\begin{Proposition}
\mylabel{real inclusion}
Let $\CCOR \subset \CCR$ be an inclusion of $\RR$-cones and
denote by $\CCOC \subset \CCC$ the inclusion of the
corresponding  complexified cones.
Let $\Delta_\RR = \diam_{\CCR} (\CCOR)^* \in [0,+\infty]$
and $\Delta_\CC = \diam_{\CCC} (\CCOC)^* \in [0,+\infty]$
 be the projective diameters
of the respective inclusions. 
Then $\Delta_\RR$ is finite iff $\Delta_\CC$ is finite.
\end{Proposition}

Proof: From the embedding in  Theorem \ref{thm embedding} 
we see that  $\Delta_\RR \leq \Delta_\CC$ which implies one direction.

To see the converse,
suppose that $\eta=\tanh \Delta_\RR/4 < 1$ and
 let $u_1,u_2\in\CCOC^*$. Possibly after rotating the polarization
 of $u_1$ and $u_2$, we may by 
 Lemma \ref{complex polarization}
assume that $u_1=x_1+i y_1$ and $u_2=x_2+i y_2$ with
\begin{equation}
 \mylabel{y bound}
 | \la \ell, y_1 \ra | \leq \eta \la \ell, x_1 \ra \ \ \
  \mbox{and} \ \ \
| \la \ell, y_2 \ra | \leq \eta \la \ell, x_2 \ra
\end{equation}
for all $\ell\in \CCPR$. Eventually applying a real rescaling of e.g.\ $u_1$
	(which does not change its polarization),
 we may by Lemma \ref{real rescaling}
assure that in addition~:
\begin{equation}
\mylabel{x bound}
 |\la \ell, x_1-x_2\ra | \leq \eta \la\ell,x_1+x_2\ra.
 \end{equation}

Let us write $u_\lambda = (1+\lambda)u_1 + (1-\lambda)u_2$, 
$\lambda\in\CC$ and similarly for $x_\lambda$ and $y_\lambda$.
In order to prove our claim it suffices to
 find a fixed open neighborhood
$U(\eta)$ of the segment $[-1;1]\subset \CC$,
depending on $\eta$ but not upon $u_1$ and $u_2$, such that
$u_\lambda\in\CCC$ for every $\lambda\in U(\eta)$.
Let $-1\leq t\leq 1$. Then
$|\la \ell,y_t\ra | \leq \eta 
       \la \ell,x_t\ra $ and we get
(with  $\ell_1,\ell_2\in\CCPR$)
       (a)~:
\[ \Re \la \ell_1,u_t\ra \la \ell_2,\baru_t \ra  =
        \la \ell_1,x_t\ra \la \ell_2,x_t \ra  +
        \la \ell_1,y_t\ra \la \ell_2,y_t \ra \geq
        (1-\eta^2)
	   \la \ell_1,x_t\ra \la \ell_2,x_t \ra\]
as well as
(b)~: $|\la \ell,u_t\ra | \leq  
       \sqrt{1+\eta^2} \la \ell,x_t\ra $.
 We also obtain
 the estimates (c)~: $|\la\ell,u_1-u_2\ra| \leq \eta \sqrt{2} \,
	   \la \ell,x_1+x_2\ra$ and (d)~:
$|\la \ell,x_t\ra | \geq 
       (1-\eta|t|) \la \ell,x_1+x_2\ra $.
Let us write $\lambda=t+z$ with $-1\leq t\leq 1$ and $z\in \CC$.
Using the expansion
  $\la \ell,u_\lambda\ra = 
   \la \ell,u_t\ra  + z \la \ell,u_1-u_2\ra$ 
   and inserting the estimates (a)-(d) we  obtain
\begin{eqnarray*}
\lefteqn{
  \Re \la \ell_1,u_\lambda\ra \la \ell_2,\overline{u_\lambda} \ra } \\[2mm]
    &\geq &
        (1-\eta^2)  \la \ell_1, x_t\ra
          \la \ell_2,x_t \ra -\\
       & &|z| \eta \sqrt{2} \sqrt{1+\eta^2} 
	     \left(
	     \la \ell_1,x_t \ra \la \ell_2,x_1+x_2\ra +
	     \la \ell_1,x_1+x_2 \ra \la \ell_2,x_t\ra\right) \\
	     & & -|z|^2 2\eta^2 
	     \la \ell_1,x_1+x_2 \ra \la \ell_2,x_1+x_2\ra\\[2mm]
	    &\geq& 
          \la \ell_1, x_t\ra \la \ell_2,x_t \ra  
    \left(2 - \left(\sqrt{1+\eta^2} + |z|
             \frac{\eta \sqrt{2} }{1-\eta |t|} \right)^2 \right).
\end{eqnarray*}
This remains positive when $|z|$ is sufficiently small
(recall that $1-\eta|t|\geq 1-\eta>0$).
The set of such $\lambda=t+z$-values thus defines an open connected
neighborhood, $U(\eta)\subset \CC$ of the segment
$[-1;1]\subset \CC$. 
Since enlarging a domain decreases hyperbolic distances, we conclude that
$\Delta_\CC\leq d_{U(\eta)} (1,-1) <\infty$. 
\Halmos\\
  
\begin{Remark}
\mylabel{C contains R}
Suppose that $T$ is a real bounded linear operator, that
$\CCR$ is a $T^{n_0}$-inner regular (some $n_0\geq 0$),
 norm-directed cone and that
the real operator $T$ maps $\CCR$ into a subcone of finite 
projective diameter in $\CCR$. 
By Birkhoff's Theorem, \cite{Bir67}, the operator
has a spectral gap. In view of Remark \ref{Remark norm-directed}
and the properties of the canonical complexification shown above,
the same conclusion follows when considering the complexified
operator acting on the canonically complexified cone. 
Our complex cone contraction thus contains the real contraction
as a special case (but, of course, with a more complicated proof).
\end{Remark}

\section{Dominated complex cone-contractions}
  \mylabel{C dominated}
  A real operator $P$ which contracts a real cone $\CCR$
contracts {\em a fortiori}
the corresponding complexified cone $\CCC$ (easy). It is then natural
to ask if this complex contraction may be preserved when
adding an imaginary part to the operator.
Many of our applications below
are cast over this idea and has lead us to
state an abstract assumption for the action upon $\CCR$ and a corresponding 
complex contraction Theorem for complexified cones~:

\begin{Assumption} 
 \mylabel{main assumption}
 Let $P\in L(X_\RR)$ be a contraction of an $\RR$-cone $\CCR$.
 Let $M\in L(X_\CC)$ be an operator acting upon the corresponding
 complex Banach space.  We say  that $M$
 is dominated by  $P$ 
  with constants
 $0\leq \gamma < \alpha\leq \beta<+\infty$ 
  provided  that
  for all $\ell,\ell_1,\ell_2\in\CCPR$ and $x,x_1,x_2\in\CCR$~:
   \begin{eqnarray}
          \Re \la \ell_1,Mx\ra  \la \ell_2,\barM x \ra &\geq &
                  \alpha \la \ell_1,Px\ra  \la \ell_2,P x \ra
		  \mylabel{domin a}\\[2mm]
   \Re \la \ell,Mx\ra  \la \ell,\barM x \ra &\leq &
                  \beta \la \ell,Px\ra  \la \ell,P x \ra
		  \mylabel{domin b}\\[2mm]
    | \Im \la \ell_1,Mx_1\ra  \la \ell_2,\barM x_2 \ra | &\leq &
                  \gamma \la \ell_1,Px_1\ra  \la \ell_2,P x_2\ra
		  \mylabel{domin c}
   \end{eqnarray}\\[-3mm]
We say that $\CCR$ is $M^{n_0}$-inner regular ($n_0\geq 0$) if there are
$x_0\in \CCR$ and $C<\infty$ so that
 \begin{equation}
       | \la \ell, M^{n_0} u \ra | \leq C  \la \ell,x_0 \ra \ \|u\|,
       \ \ \ \forall u\in X_\CC , \ \ell\in \CCPR.
             \mylabel{aux condition}
 \end{equation} \mbox{}
\end{Assumption}

\begin{Remark}
\mylabel{remark generating}
The above conditions 
are $\RR_+$-invariant and also
stable when taking convex combinations. It thus suffices
to verify that these conditions hold for  subsets,
$V\subset \CCR$ and  $W\subset \CCPR$ 
which are generating for the cone and the dual cone, respectively, i.e.\ for which~:
 \[ \CCR=\Cl \Conv (\RR_+\times V) =
   \left\{ x \in X_\RR : \la \ell,x\ra \geq 0,\ 
      \ \forall \; \ell \in W\right\}.\]
 Note also that when $\CCR$ is inner regular 
  the last condition (\ref{aux condition})
  is automatically satisfied
 (choose for $x_0$ an interior point in $\CCR$).
\end{Remark}

When $\gamma=0$ an operator $M$ verifying
the above assumption is essentially
real. Possibly after multiplication with a complex
phase the operator maps $\CCR$ into $\CCR$ itself.
The above condition then reduces to the real
cone-dominated condition of Theorem \ref{real domin gap}.
Our goal is here to show that the conclusion of that Theorem
also applies when $M$ is allowed to have a non-trivial 
imaginary part. 
It turns out that the allowed `amount' of
imaginary part depends on the rate of contraction of $P$.

\begin{Theorem}
\mylabel{C domin sp gap}
Let $\CCR\subset X_\RR$ be a proper convex cone and
let $P:\CCSR\rr\CCSR$ be a strict cone-contraction, i.e.\ 
$\Delta_P = \diam_{\CCR} P (\CCSR) < +\infty$.
Write $\CCC$ for the canonical complexification of $\CCR$.
Suppose that $M\in L(X_\CC)$ is $P$-dominated 
(Assumption \ref{main assumption}) with constants that satisfy~:
  \begin{equation}
  \gamma \; \cosh \frac{\Delta_P}{2} < \alpha.
  \label{main condition}
  \end{equation}
Then $M: \CCSC \rr \CCSC$ and $\diam_{\CCC} M(\CCSC) <+\infty$.
If, in addition, $\CCR$ is $M^{n_0}$-inner regular ($n_0\geq 0$) 
and of uniformly bounded sectional aperture then
 $M$ has a spectral gap.
\end{Theorem}
Proof: 
Let $u\in\CCSC$ and $\ell_1,\ell_2\in \CCPR$. We
write $\eta=\tanh\frac{\Delta_P}{4}<+\infty$.
The first step is to establish the following inequality
(which, in particular, implies that $M:\CCSC\rr \CCSC$)~:
\begin{equation}
   \mylabel{eq ell u}
   \Re \la \ell_1,M u\ra \la \ell_2,\overline{Mu} \ra \geq
     \left( \frac{\alpha}{\cosh ({\Delta_P}/{2})} - \gamma\right) 
             |\la \ell_1,Pu\ra | 
             |\la \ell_2,Pu\ra |  .
\end{equation}
We will use polarization twice to achieve this.
  First, write $u=e^{i\theta}(x+iy)$ with
   $\theta\in\RR$ and $x\pm y\in \CCSR$.  Then
 \begin{eqnarray*}
     \lefteqn{\la \ell_1,Mu\ra   \overline{\la \ell_2,Mu\ra}} \\
       &= &\la \ell_1,M( x\!+\!iy)\ra   \la \ell_2,\barM(x\!-\!iy)\ra =\\
       &= &\left[\la \ell_1,M x\ra   \la \ell_2,\barM x\ra 
           + \la \ell_1,M y\ra   \la \ell_2,\barM y\ra \right]
        + i \left[\la \ell_1,M y\ra   \la \ell_2,\barM x\ra 
           - \la \ell_1,M x\ra   \la \ell_2,\barM y\ra \right]\\
       &=&
    \frac12 \left[ \la \ell_1, M(x\!+\!y)\ra \la \ell_2, \barM(x\!+\!y)\ra
     \!+\! \la \ell_1, M(x\!-\!y)\ra \la \ell_2,\barM( x\!-\!y)\ra \right] + \\
  &&  \frac{i}{2}\left[ \la \ell_1,M(x\!+\!y)\ra \la \ell_2,\barM(x\!-\!y)\ra
     \!-\! \la \ell_1, M(x\!-\!y)\ra   \la \ell_2,\barM(x\!+\!y)\ra \right] \\
     &\equiv& \frac12 [A] + \frac{i}2 [B].
 \end{eqnarray*}
Since $x\pm y\in\CCR$ we may use  
inequality (\ref{domin a}) of our assumption
to deduce~:
\begin{eqnarray*}
     \;\Re\; A 
    &\geq&\ 
          \alpha \; \la \ell_1,P(x+y) \ra
          \la \ell_2,P(x+y) \ra +
	  \alpha \; \la \ell_1,P(x-y) \ra
          \la \ell_2,P(x-y)\ra\\
    &=& 
          2 \alpha \; \la \ell_1,Px \ra
          \la \ell_2,Px \ra +
	  \alpha \la \ell_1,Py \ra
          \la \ell_2,Py\ra\\
    &=& 
          2 \alpha \; \Re  \;
	   \la \ell_1,P(x+iy) \ra
          \la \ell_2,P(x-iy) \ra \\
    &=& 
          2 \alpha  \;\Re  \;
	   \la \ell_1,Pu \ra
          \la \ell_2,P\baru \ra .
\end{eqnarray*}
For the second term we have by (\ref{domin c})
\begin{eqnarray*}
     |\; \Im \;B| 
    &\leq&\ 
       \gamma \;  \la \ell_1,P(x+y) \ra
          \la \ell_2,P(x-y) \ra +
	  \gamma\;\la \ell_1,P(x-y) \ra
          \la \ell_2,P(x+y)\ra\\
    &=& 
          2 \gamma \;(\la \ell_1,Px \ra
          \la \ell_2,Px \ra -
	  \la \ell_1,Py \ra
          \la \ell_2,Py\ra)\\
    &\leq & 
        2 \gamma \;
	   |\la \ell_1,P(x+iy) \ra|\;
          |\la \ell_2,P(x-iy) \ra| \\
    &=& 
          2 \gamma \;
	   |\la \ell_1,Pu \ra|\;
          |\la \ell_2,Pu \ra| ,
\end{eqnarray*}
where for the last inequality we used Schwarz' inequality.
   From these two estimates we get~:
\begin{equation}
    \mylabel{Mu Mu}
    \Re \la \ell_1,Mu\ra   \overline{\la \ell_2,Mu\ra}\geq
          \alpha \;\Re 
	   \la \ell_1,Pu \ra
          \la \ell_2,P\baru \ra -
          \gamma \;
	   |\la \ell_1,Pu \ra|\;
          |\la \ell_2,Pu \ra| .
\end{equation}

We note that (\ref{Mu Mu}) is here independent of the choice of
polarization.
Since  $x\pm y\in\CCSR$ we see that
the elements $P(x+y)\in \CCSR$ and $P(x-y)\in \CCSR$ are at
a projective distance not exceeding $\Delta_P$. We may then use
Lemma \ref{complex polarization} to rotate the polarization again
and write $Pu=e^{i\alpha}(x'+iy')$ where
$|\la \ell,y'\ra| \leq \eta \la \ell,x'\ra$
for all $\ell\in\CCPR$.
But then
\[ \Re \la \ell_1,Pu \ra \la \ell_2,P\baru \ra \geq
          \la \ell_1,x' \ra \la \ell_2,x' \ra +
	  \la \ell_1,y' \ra \la \ell_2,y'\ra
	  \geq (1-\eta^2) 
          \la \ell_1,x' \ra \la \ell_2,x' \ra .\]
We also obtain 
$|\la \ell,Pu\ra|
   = \sqrt{
          \la \ell,x' \ra^2  +
	  \la \ell,y' \ra^2} \leq
	  \sqrt{1+\eta^2}\;
	  \la \ell,x'\ra$ so that
\[ \Re \la \ell_1,Pu \ra \la \ell_2,P\baru \ra \geq
	   \frac{1-\eta^2}{1+\eta^2}
	    |\la \ell_1,Pu \ra|\; |\la \ell_2,Pu \ra| =
     \left(\cosh \frac{\Delta_P}{2}\right)^{-1} 
	   |\la \ell_1,Pu \ra|\; |\la \ell_2,Pu \ra| .
   \]
Together with (\ref{Mu Mu}) this establishes (\ref{eq ell u}).

In order to obtain an estimate for 
$\diam_{\CCC} M(\CCSC)$
we also need the following inequality~:
\begin{equation}
   |\la \ell, Mu\ra | \leq \sqrt{\beta + \gamma} \; \
   |\la \ell, Pu\ra|, \ \ \forall \ell\in \CCPR, u\in \CCC.
   \mylabel{Gamma ineq}
   \end{equation}
This follows by setting $\ell_1=\ell_2=\ell$ in the expression for $A$
and $B$ above 
 and using the upper bounds (\ref{domin b}) and (\ref{domin c})
 of our Assumption~:
\[ A \leq
   \beta (
       \la \ell,P(x+y)\ra^2 +
       \la \ell,P(x-y)\ra^2)
   = 2\beta (
       \la \ell,Px\ra^2 +
       \la \ell,Py\ra^2) = 2\beta \;
       |\la \ell, Pu\ra|^2
\]
and the bound
$ |B|=|\Im B| \leq 2\gamma \; |\la \ell, Pu\ra|^2$ as before.\\

Consider  $u_1,u_2\in\CCSC$. 
Using the polarization identity,
 Proposition \ref{equiv C},
  we may assume that
 $u_1=x_1+iy_1$  with $x_1\pm y_1\in \CCSR$ so that
$|\la \ell, P y_1|\leq \la \ell,Px_1\ra$.
Then also
$\la \ell,P x_1\ra \leq |\la \ell, P u_1 \ra| \leq \sqrt{2} \la \ell,Px_1\ra$
and with the same bounds  for $u_2=x_2+iy_2$.
Through a real rescaling,
Lemma \ref{real rescaling},
we may also assume that
$| \la \ell,P (x_1-x_2)\ra | \leq \eta \la \ell,P (x_1+x_2)\ra $.
We also write $u_\lambda=(1+\lambda)u_1+(1-\lambda) u_2$ with 
$\lambda=t+z$, $-1\leq t \leq 1$ (and similarly
for $x_\lambda$ and $y_\lambda$). 
	By the choice of polarization 
$x_t\pm y_t\in\CCR$ so that $u_t\in \CCC$, i.e.\ belongs to the complex cone
 for all $-1\leq t \leq 1$.
We want to show that when $|z|$ is small enough the same is true for
$Mu_{t+z}$.

 First note that
  $\la \ell,Px_t\ra \geq (1-\eta|t|) \; \la \ell,P(x_1+x_2)\ra$.
Applying the inequality (\ref{Gamma ineq}) we 
deduce that
$\ds |\la \ell,M u_t \ra | \leq \sqrt{2\;(\beta+\gamma)} \;
      \la \ell,P x_t\ra$ and 
$\ds |\la \ell, M (u_1-u_2)\ra | \leq 
    \sqrt{2\;(\beta+\gamma)} \la \ell, P(x_1+x_2)\ra \leq \sqrt{2\;(\beta+\gamma)} \;
    \frac{\la \ell, P x_t\ra}{1-\eta|t|}$.
Using (\ref{eq ell u}) on $u_t$ and the expansion
$\ds \la \ell,Mu_{t+z}\ra = \la \ell, M u_t\ra + z \la \ell,M(u_1-u_2)\ra$,
%proceed as in Proposition \ref{real inclusion} and
we obtain the inequality
\begin{eqnarray}
  \lefteqn{\Re\; \la \ell_1, M u_{t+z} \ra\; 
     \la \ell_2, \overline{M u_{t+z}}\ra \;} \\[3mm]
        & \geq  & \Re\,
  \la \ell_1, Mu_t\ra\, 
  \la \ell_2, \overline {Mu_t}\ra\  \nonumber \\
  &  & - |z| \  \left(\ \strut
  |\la \ell_1,M u_t\ra | \ |\la \ell_2,M (u_1-u_2)\ra | \; +\;
  |\la \ell_2,M u_t\ra | \ |\la \ell_1,M (u_1-u_2)\ra |  \
          \right)
   \nonumber \\
  &  & - |z|^2 \ 
  |\la \ell_1,M (u_1-u_2)\ra | \ 
  |\la \ell_1,M (u_1-u_2)\ra |   \nonumber \\[2mm]
   & \geq &\la \ell_1, P x_t\ra \, \la \ell_2 P x_t \ra \ \times
    \left(\frac{\alpha}{\cosh({\Delta_P}/{2})}-\gamma 
    + 2 (\beta+\gamma)  - 2(\beta+\gamma)  \left(
    1 + \frac{|z|}{1-\eta|t|} \right)^2
        \right) . \nonumber
\mylabel{C diam}
\end{eqnarray}
The set of $t+z$ values, $-1\leq t\leq 1$
 for which the latter quantity is non-negative contains an open neighborhood,
 $U=U(\alpha,\beta,\gamma,\Delta_P)$
 of the segment $[-1;1]\subset \CC$. It follows that
 \[ \diam_{\CCC} M(\CCSC) \leq d_U(-1,1) <+\infty.\]

For the last assertion note that the condition (\ref{aux condition})
implies that 
\[ \Re 
\la \ell_1, x_0+M^{n_0} u\ra 
    \la \ell_2, x_0+\overline{M^{n_0} u}\ra  \geq
    \la \ell_1,x_0\ra     \la \ell_2,x_0\ra
       (1 - 2 C \|u\| -  C^2 \|u\|^2)
       \] which remains
   non-negative when $\|u\| \leq \frac{1}{C}(\sqrt{2}\, - 1)$.
Thus $\CCC$ is $M^{n_0}$-inner regular. By 
 Proposition \ref{C regular cone}, $\CCC$ is  also 
 of uniformly bounded sectional aperture
so the conclusion follows from our spectral gap theorem, 
Theorem \ref{spectral gap}.\Halmos\\

\section{Applications}
\label{sec Applications}
The most striking application is also the simplest. A complex
Perron-Frobenius Theorem~:
 
 \begin{Theorem}  
 Let $A\in M_n(\CC)$ and suppose
there is $0<c<+\infty$ for which 
$ |\im\, A_{ij}\overline{A}_{mn}|
    <  c \leq  \re\, A_{ij}\overline{A}_{mn}$ for all indices. Then $A$ has
    a spectral gap.
  \mylabel{Thm PF complex}    
\end{Theorem}
    Proof~:  The cone $\CCR=\RR_+^n$ is regular in $\RR^n$.
By Corollary \ref{Corollary Cplus}, $\CCC=\CC_+^n$ is regular in $\RR^n$.
    We will compare
    $M$ with the constant matrix
     $P_{ij}\equiv 1$ with respect to the real cone $\RR_+^n$.
     As in Example \ref{ex real PF}(1), $\Delta_P=0$.
    The canonical basis and its dual generates the cone and its dual,
    respectively, cf.\ remark \ref{remark generating}.
    The constants from Assumption \ref{main assumption} then become
    (sups and infs over all indices)
    (a) $\alpha=\inf \Re A_{ij}\barA_{kl}$, 
        (b) $\beta=\sup \Re A_{ij}\barA_{kl}$
	and (c) $\gamma=\sup |\Im A_{ij} \barA_{kl}|$. Our spectral
	gap condition of Theorem \ref{C domin sp gap} simply reads $\gamma<\alpha$ and
    by finiteness of $n$ this is equivalent to the stated assumptions on $A$.
    \Halmos\\

In the following, denote by 
$\osc(h)=\ess \sup(h) - \ess \inf (h)$ the essential oscillation of
a real valued function $h$ on a measured space.
Theorem \ref{Thm PF complex} may (almost) be viewed as a special case of
the following complex version of a result of Jentzsch
\cite{Jen12}~:

\begin{Theorem} \mylabel{Complex Jentzsch}
Let $(\Omega,\mu)$ be a measure space and
let $X=L^{p}(\Omega,\mu)$, with $1\leq p\leq +\infty$.
Let  $h\in L^p$, $h>0$ a.e.\ and
$m\in L^q$, $m>0$ a.e.\ with $\frac{1}{p}+\frac{1}{q}=1$
 so that $0<\int h\, m\, d\mu< +\infty$,
cf. Example \ref{ex real Jentzsch}(2).
Given $g\in L^\infty(\Omega\times \Omega)$
we define the integral operator,
$M_g\in L(X)$~:
  \begin{equation}
    M_g \phi(x)=h(x)\int_\Omega e^{g(x,y)} \phi(y)\, m(y) \,\mu(dy).
  \end{equation}
Set $\theta=\osc (\Im\; g)$ and $\Lambda=\osc(\Re\; g)$.
 Suppose that $\theta<\pi/4$ and that 
$\;\tan \theta < \exp(-2 \Lambda)$. Then $M_g$ has a spectral gap.
\mylabel{Jentzsch complex}
\end{Theorem}

Proof: As in Example \ref{ex real Jentzsch}(2) 
we consider the $\RR$-cone
$\CCR=\{\phi\in X_\RR :  \phi \geq 0 \ {\rm (a.e.)} \}$ and
%of functionals $\ell_x\in\CCPR$, $\ell_x(\phi)=\phi(x)$.
we compare with $ P\phi=h \int_\Omega \phi \, m\,d\mu$.
We have that $P:\CCSR\rr \CCSR$ and $\Delta_P=0$.
We obtain the following estimate for  the constants
\[
   \Re \;e^{g(x,y) + \overline{g(x',y')}} \geq
     \alpha \equiv e^{2 \ess \inf \Re\; g} \cos\theta \ \ \mbox{and} \ \ 
\Im \;e^{g(x,y) + \overline{g(x',y')}} \leq
     \gamma \equiv e^{2 \ess \sup \Re\; g} \sin\theta .\]
 Since $|\la \ell, M_g \phi\ra | \leq
\la \ell,h\ra e^{2 \ess \sup \Re g } \|m\|_q \|\phi\|_p$,
$\phi\in X_\CC$ and $\ell\in \CCPR$,
  the cone $\CCR$ is $M_g$-inner regular.
 As shown in Example \ref{ex real Jentzsch} (2) the real cone
 has bounded sectional aperture so by Proposition
\ref{C regular cone} (3), the same is true for the complexified cone.
 The spectral gap condition of Theorem \ref{C domin sp gap}
     then  translates into the stated condition on
 $\theta$ and $\Lambda$. \Halmos\\

\section{A complex Kre\u{\i}n-Rutman Theorem}
%\mylabel{section C finite}
Let $X$ be a complex Banach space. We denote by 
$\grass(X)$ denote the set of
complex planes in $X$, i.e.\ subsets of the form $\CC x+\CC y$ with 
$x$ and $y$ independent vectors in $X$.  If we write $S(X)$ for the unit
sphere in $X$ then
\[ d_2(F,F')=\dist_H (F\cap S(X),F'\cap S(X)), \ \ F,F'\in \grass(X)\]
defines a metric on  $\grass(X)$.
In the following let us fix a norm on  $\CC^n$.
 The choice may affect the constants below but is otherwise immaterial.
The space $(\grass(\CC^n),d_2)$ is then
a sequentially compact metric space.

\begin{Lemma} Let $V\subset X$ be a $\CC$-cone and let $F\in \grass(X)$.
Suppose there is $u\in F$, $r>0$ such that $B(u,r) \cap V=\emptyset$.
Then $V_F=F\cap V$ has at most  $1+\frac{\|u\|}{r}$ bounded aperture.
\end{Lemma}

Proof: Let $m\in (F)'$ be a linear functional
with $u\in \ker m$ and $\|m\|=1$. Choose
$x\in F$ for which $|\la m,x\ra|=\|x\|$. 
 If $ax+bu\in V_F$ then
$u+\frac{a}{b} x\notin B(u,r)$ so that $|b|\leq \frac{|a|}{r} \|x\|$
and therefore,
\[ \| ax+bu\| \leq |a| \|x\| (1+\frac{\|u\|}{r}) = 
      |\la m,a x+b u\ra| (1+ \frac{\|u\|}{r}). \]
The 2-dimensional space $F$ is
 spanned by $u$ and $x$ so $K(V_F) \leq 1+\frac{\|u\|}{r}$.\Halmos\\

\begin{Lemma}
\mylabel{unif K aperture}
Let $V\subset \CC^n$ be a $\CC$-cone. Then there is $K<\infty$ so that
$V$ is of $K$-bounded sectional aperture.
\end{Lemma}
Proof: Suppose that this is not the case. Then we may find a sequence $F_n$ of
planes for which the aperture $K(V\cap F_n)$ diverges. Taking a subsequence
we may assume that $F_n$ converges in 
  $\grass(X)$
 to a plane $F$. As $V$ is proper, $V\cap F$
is a strict subset of $F$. Thus there is $u\in F-V$. But $V$ is closed in $\CC^n$
 so there is
$r>0$ so that $B(u,r)$ is disjoint from $V$ as well.
 Given another complex plane, $F'$, we may find $u'\in F'$ for which
 $\|u-u'\| \leq \|u\| d_2(F,F')$. When $F$ and $F'$ are close enough,
 $r'=r-\|u\| d_2(F,F')>0$ and 
  $B(u',r')$ is also disjoint from $V$. By our previous
  Lemma, $V\cap F'$ is of aperture
  not exceeding $1+{\|u\|}/{(r-\|u\| d_2(F,F'))}$. But this contradicts the
  divergence of $K(V \cap F_n)$ as $F_n\rr F$.\Halmos\\
 
\begin{Lemma} 
\mylabel{Wv cone}
Let $V\subset \CC^n$ be a  $\CC$-cone and let $W\subset V$
be a closed complex sub-cone with $W^* \subset \Int V$.
Then there is $\Delta=\Delta(W,V) < +\infty$ such that for $x,y\in W^*$~:
  \[ d_W(x,y) < +\infty \Rightarrow d_V(x,y) \leq \Delta.\]
\end{Lemma}
Proof: 
%By the previous Lemma, $V$ is of $K$ bounded sectional aperture for some
%$K<\infty$. 
We denote by $\pi: \CC^n-\{0\} \rr \CC P^{n-1}$ the canonical 
projection to complex projective space. We consider
 $\CC P^{n-1}$ as a metric space with
 the metric $d_{\CC P^{n-1}}$ as in (\ref{proj metric}).
The projected image, $\pi(W^*)$,
 is  compact in the open set $\pi(\Int V^*)\subset \CC P^{n-1}$
so  there is  $\epsilon=\epsilon(W,V)>0$ for which the 
$\epsilon$-neighborhood of $\pi(W^*)$   is
contained in $\pi(V^*)$.

%Then the $\epsilon/3$-neighborhoods of $G$ covers $\CC P^{n-1}$.
Let $x,y\in W^*$ be linearly independent
and suppose that $d_{W}(x,y)<+\infty$.
Let $F\in \grass(\CC^n)$ be the  
complex plane containing $x$ and $y$.  Denote by
$C$ the connected set in $\pi(F^*)$ containing $x$ and $y$.
Let $\xi_i\in C$, $i\in J$ be an
$\epsilon/3$-maximally separated set in $C$.
Thus, the balls
 $B(\xi_i,\frac{\epsilon}{6})$, $i\in J$ are all disjoint and
$\bigcup_{i\in J}B(\xi_i,\frac{\epsilon}{3})=\CC P^{n-1}$. 
The cardinality of $J$ is bounded by a constant depending on
$\epsilon$ only.
Then $B(\xi_i,\frac{2\epsilon}{3}) \subset \pi(V^*)$, $i\in J$ so 
by Lemma \ref{XP lower bound}(2) each 
$B(x_i,\frac{\epsilon}{3})$, $i\in J$
 is of radius  not greater than $\log \frac{1+1/2}{1-1/2}=\log 3$
  for the $d_V$-metric.
Also $\cup_{i\in J} B(\xi,\frac{\epsilon}{3})$  contains $C$ 
which is connected.
It follows that
$d_V(x,y)$ does not exceed 
   $2 \, \log 3\; \Card (J)$.
which is bounded by a constant depending
on $\epsilon$ only.\Halmos\\

\begin{Theorem}
 \mylabel{Thm finite dim}
 Let $V\subset \CC^n$ be a closed subset which is 
$\CC$-invariant and contains
no complex planes (in terms of Definition \ref{def complex cone}, $V$ is
  a  $\CC$-cone).
Suppose that $A:\CC^n\rr \CC^n$ is a linear map for 
which $A (V^*) \subset \Int V$.
Then $A$ has a spectral gap. 
\end{Theorem}
Proof~: We write $W=A (V)$ for the image of $V$ and use the 
notation and constants from
the two previous Lemmas.
First note that for $x,y\in V$,
\[ d_V(x,y)<\infty \Rightarrow d_V(Ax,Ay) \leq \eta\; d_V(x,y) \ 
  \ \ \mbox{and}  \ \ \  d_V(Ax,Ay) \leq \Delta.\]
To see this note that when $d_V(x,y)<\infty$ then $d_V(Ax,Ay)<\infty$ so
by Lemma \ref{Wv cone},
$d_V(Ax,Ay)\leq \Delta$.
If $\calC^*$ denotes the connected component of $F\cap V^*$
 containing $x$ and $y$
then also $\diam_V A(\calC^*) \leq \Delta$.
By Lemma \ref{lemma contraction}, $d_V(Ax,Ay) \leq \eta\, d_{\calC}(x,y)
 \leq \eta d_V(x,y)$. Iterating this argument we see that
 $\diam_V A^n (\calC^*) \leq \Delta \eta^{n-1}$, $n\geq 1$.
By Lemma \ref{unif K aperture},
 $V$ is of $K$-bounded sectional aperture, so Lemma \ref{XP lower bound}(1)
assures that
 $\diam_{\CC P^{n-1}} A^n (\calC^*)\leq 2 K \Delta \eta^{n-1}$.
 Fix $n_1<+\infty$ so that
 $2K \Delta \eta^{n_1-1}\leq \epsilon/3$.
 
 Now let $\xi_i$, $i\in J$ be an $\epsilon/3$-maximally separated set in $W$.
 Setting $V_i= \pi^{-1} B(\xi_i,\epsilon)$ with $i\in J$ we see that
 $\diam_{\CC P^{n-1}} A^n V_i^* \leq \epsilon/3$, $n\geq n_1$.
 It follows that there is a map, $\tau:J\rr J$ so that
 $A^n V^*_i \subset W_{\tau(i)}\equiv \pi^{-1}
  B(\xi_{\tau(i)},2\epsilon/3)$, $n\geq n_1$.
 Since $J$ is of finite cardinality, $\tau$ must have a cycle.
 Thus, there are $i_1\in J$ and $n_1<+\infty$ for
 which $A^{n_1} (V_{i_1}) \subset W_{i_1}$. The cone $W_{i_1}$
 is regular (easy) and of bounded 
 diameter in $V_{i_1}$ 
  so $A^{n_1}$ has a spectral gap and therefore
  also $A$.\Halmos\\

When the operator is sufficiently regular one may weaken 
the assumptions on the contraction and the
outer regularity of the cone. This is illustrated
by the following complex version of   a
theorem of Kre\u{\i}n and Rutman \cite[Theorem 6.3]{KR50}~:
\begin{Theorem}
\mylabel{complex Krein Rutman}
 Let $\calC\subset X_\CC$ be a $\CC$-cone in the
 Banach space $X_\CC$. Let $A\in L(X_\CC)$ be a  quasi-compact
 operator or a compact operator of strictly positive spectral radius and suppose that
 $A$ verifies
    \begin{equation}
      A : \calC^* \rr \calC^o,
    \end{equation}
 Then $A$ has a spectral gap.
\end{Theorem}
Proof~:
Let $P$ be the spectral projection associated with eigenvalues
on the spectral radius circle, $\{\lambda\in \CC: |\lambda| = r_{\rm sp}(A)\}$.
By hypothesis  ${\rm im} P$ is finite dimensional and we may
find $\theta\in\RR$ such that
  \[ r_{\rm sp} (A(1-P)) < \theta < r_{\rm sp}(A).\]
We claim that $\CCS \cap {\rm im} P$ is non-empty~:
Let $x\in \calC^*$ and  
define $e_n=A^n x/\|A^nx\|\in \CCS, \ n\in\NN$.

Suppose first that $Px\neq 0$. Then 
$\lim_{n\rr\infty} {\|A^n(1-P)x\|}/{\|A^nPx\|} =0$ so that the distance
between $e_n$ and ${\rm im} P$ tends to zero. Since $\; {\rm im}\; P\; $ is locally compact
and $e_n$ is bounded we may extract
a convergent sub-sequence $e^*=\lim e_{n_k} \in {\rm im} P \cap \CCS$.
 Suppose instead that $Px=0$ then $Ax\in \calC^o$
so there is $r>0$ for which $B(Ax,r)\in \calC$. We may then replace $x$ by
$Ax+u$ where $u\in {\rm im} P$, $\|u\|<r$ and we are back in the first case.
Thus $\calC_P^*=\calC^* \cap {\rm im} P \neq \emptyset$.
 Now,
  \[ A : \calC_P^* \rr (A \;\calC^*)\cap {\rm im} P \subset 
                 \calC^o \cap {\rm im} P = \calC_P^o ,\]
the latter for the topology in ${\rm im} P$. In particular,
$\calC_P^o$ is non-empty so
 $\calC_P$ is an inner regular $\CC$-cone in a finite
dimensional space and $A: \calC_P^* \rr \calC_P^o$. We may then apply
the finite dimensional contraction theorem, Theorem \ref{Thm finite dim},
to $A_P=A_{|{\rm im P}} \in L({\rm im} P)$. It follows that $A_P$, whence also
our original operator $A$ has a spectral gap.\Halmos\\

 \begin{Remark}
 \mylabel{KR remark}
 In the real cone version (replacing $\CC$ by $\RR$)
  of theorem  \ref{complex Krein Rutman}
 it is not necessary to assume that
 the spectral radius of $A$ is strictly positive. This forms part of the
 conclusion. To see this
 pick $x\in \calC^*$ of norm one. Then $Ax\in \calC^o$ so there is
$\lambda>0$ for which $B(Ax,\lambda)\subset \calC$.
Therefore, $Ax-\lambda x\in \calC$ and then also
$B(A^2 x,\lambda^2) = A(Ax-\lambda x) + \lambda B(Ax,\lambda) \subset \calC$
by the  properties of an $\RR$-cone. More generally,
$B(A^n x, \lambda^n) \subset \calC$.
As $0\in \partial \calC$ it follows that
  \[ r_{\rm sp} (A) \geq \limsup \sqrt[n]{|A^nx|} \geq \lambda >0 . \]
The fact that this conclusion is non-trivial is illustrated e.g.\ by
the operator, $A\phi(t)=\int_0^s \phi(s)\; ds$, $0\leq t \leq 1$,
which is compact when acting upon
$\phi\in X=C^0([0,1])$. It contracts (though  not strictly)
the cone of positive elements but has spectral radius zero.

In the complex setup, if one assumes that $\calC$ is of $K$-bounded sectional
aperture then strict positivity of $r_{\rm sp}(A)$ 
also comes for free~:
Suppose that  $x\in\calC$, $|x|=1$
and $B(Ax,r)\subset \calC$, $r>0$. Then
$Ax+\lambda x\in \calC^*$, $\forall |\lambda|<r$ and also
$A^{n+1}x+\lambda A^n x\in \calC^*$ for such $\lambda$-values.
By Lemma \ref{lemma inner regular} we see that
$|A^{n+1}x| \geq \frac{r}{K} |A^nx|>0$ from which
$r_{\rm sp} (A) \geq \frac{r}{K}>0$.
\end{Remark}

\section{A complex  Ruelle-Perron-Frobenius Theorem}
\label{sec RPF}
The Ruelle-Perron-Frobenius Theorem,
\cite{Rue68,Rue69,Rue78} (see also \cite{Bow75}),
ensures a spectral gap
for certain classes of real, positive
operators with applications in statistical
mechanics and dynamical systems.
Ferrero and Schmitt \cite{FS79,FS88} used Birkhoff's Theorem
on cone contraction 
to give a conceptually new proof of
the Ruelle-Perron-Frobenius Theorem.
See also \cite{Liv95} and 
\cite{Bal00} for further applications in dynamical systems.
We present here a generalization to a complex setup.

Let $(\Omega,d)$ be a metric space of finite diameter, $D<+\infty$.
When $\phi:\Omega\rr \RR$ (or $\CC$) we denote by 
$\Lip(\phi)=\sup_{x\neq y} |\phi_x-\phi_y|/d(x,y)\in [0,+\infty]$ the associated Lipschitz
constant and by $|\phi|_0$ the supremum. Then 
$X_\RR=\{\phi:\Omega\rr \RR \ | \  \|\phi\| \equiv |\phi|_0+\Lip(\phi)< +\infty\}$ 
(and similarly for $X_\CC$)
is a Banach algebra.

Let $U\subset \Omega$ and let $f:U\rr \Omega$ be an unramified
covering map of $\Omega$ which is uniformly  expanding. For simplicity, we
will take it to be of finite degree (it is an instructive exercise
to extend Theorem
\ref{Thm RPF} below
to maps of countable degree).
 More precisely, we assume
that there is $0<\rho<1$ and a finite index set $J$ 
so that for every couple $y,y'\in \Omega$
we have a pairing $\calP(y,y')=\{(x_j,x'_j): j\in J\}$
of the  pre-images,
  $f^{-1}(y)=\{x_j\}_{j\in J}$ and
   $f^{-1}(y')=\{x'_j\}_{j\in J}$, for which
$d(x_j,x'_j)\leq \rho\; d(y,y')$, $ j\in J$.\\

Fix an element $g\in X_\CC$ and define for
$\phi\in C^0(M)$  (or $\phi\in X_\CC$)~:
  \[M_g\phi  (y)=\sum_{x: f(x)=y} e^{g(x)} \phi(x), \ y\in \Omega.\]
The norm of $M_g$ when acting upon $C^0(M)$ (in the uniform norm) is given by
 \nc{\NNN}{|\!|\!|}
 \[ \NNN M_g \NNN_0=\sup_{y\in \Omega}
  \sum_{x: f(x)=y} e^{\Re g(x)},\]
and a straight-forward calculation shows that $M_g\in L(X_\CC)$ with
 $ \|M\| \leq \NNN M \NNN_0 (1+\rho\; \Lip \; g)$.\\

\begin{Theorem}
\label{Thm RPF}
Denote $a=\Lip \;\Re\; g$, $b=\Lip\; \Im\; g$ and $\theta= \osc\; \Im\; g$.
Suppose that 
\[ \left( \theta + \frac{2 \,\rho^2\, D\, b}
                   {1-\rho + \rho^2\, D\, a} \right) \
          \exp \left( 1\, + \, \rho\, \frac{1+\rho}{1-\rho} \; D\, a \right)
	  \ \frac{4}{1-\rho} < 1 .\]
Then $M_g \in L(X_\CC)$ has a spectral gap.
\end{Theorem}
 
Proof~:
We will compare $M_g$ with the real operator $P=M_{\Re\; g}$.
For $\sigma>0$ the set,
\begin{equation} 
\calC_{\sigma,\RR}= \{\phi : \Omega\rr \RR_+ \ | \
          \la \ell_{y,y'},\phi\ra  \equiv
	  \phi(y) - e^{-\sigma d(y,y')}\phi(y')\geq 0, \ \forall y,y'\in \Omega\} ,
	  \mylabel{RPF cone def}
\end{equation}
defines a proper convex cone in $X_\RR$ which in addition is regular.
Inner regularity~: Let  $\bfone(x)\equiv 1$, $x\in\Omega$
and $h\in X_\RR$.
Then $\bfone+h\in\calC_{\sigma,\RR}$  provided $\Lip\, h/(1-|h|_0)\leq \sigma$.
Whence $B(\bfone, \min(\sigma,1)) \subset \calC_{\sigma,\RR}$.
Outer regularity~:
Pick $x_0\in \Omega$ and set $\ell_0(\phi)= \phi(x_0)$.
For $\phi\in \calC_{\sigma,\RR}$ we have
$\Lip\, \phi \leq \sigma |\phi|_0$ so that
$\|\phi\| \leq (1+\sigma)|\phi|_0 \leq
      (1+\sigma)e^{\sigma D}  \ell_0(\phi)$,
and this shows outer regularity.
 
Let   $0<\sigma'<\sigma$ and  $\phi_1,\phi_2\in \calC_{\sigma',\RR}^*$.
As in (\ref{Birkhoff beta})  let 
$\beta_\sigma(\phi_1,\phi_2) = \inf \{ \lambda > 0 : 
 \lambda \phi_1 - \phi_2 \in \calC_{\sigma,\RR}\}$.
A calculation using the defining properties of the cone-family yields~:
\[ \beta_\sigma(\phi_1,\phi_2) \leq \sup_{d>0} \;
            \frac{1-\exp(-(\sigma+\sigma')d)}
	         {1-\exp(-(\sigma-\sigma')d)} \
      \sup_{y\in\Omega}\; \frac{\phi_2(y)}{\phi_1(y)}
      \leq \frac{\sigma+\sigma'}{\sigma-\sigma'}\;
            \sup_{y\in\Omega}\; \frac{\phi_2(y)}{\phi_1(y)},\]
and we get  the following bound for the diameter 
 $\Delta_\RR 
 = \diam_{\calC_{\sigma,\RR}} \calC_{\sigma',\RR}^* $,
cf.\ (\ref{Birkhoff def})~:
\[ \Delta_\RR 
    \leq 2\; \log \;\frac{\sigma+\sigma'}{\sigma-\sigma'} +
            \sup_{y,y'\in\Omega}\;\log 
	     \frac{\phi_2(y)}{\phi_1(y)}
	     \frac{\phi_1(y')}{\phi_2(y')} 
    \leq 2\; \log \;\frac{\sigma+\sigma'}{\sigma-\sigma'} +
    2\; D\; \sigma' < +\infty.\]
The injection 
$\calC_{\sigma',\RR} \hookrightarrow \calC_{\sigma,\RR}$ is thus a 
uniform contraction for the respective projective metrics.
Given $\phi\in \calC_{\sigma,\RR}$ and using the pairing $\calP(y,y')$
we get for the operator $P=M_{\Re g}$~:
\[ P\phi(y) = 
   \sum_{x:f(x)=y} e^{\Re\; g(x)} \phi(x) \geq
   \sum_{x':f(x')=y'} e^{\Re\; g(x') -(a +\sigma) d(x,x')} \phi(x') \geq
    e^{-\rho(a+\sigma) d(y,y')} P\phi(y').\]
 This implies that 
  $ P : \calC_{\sigma,\RR}\rr
\calC_{\sigma',\RR}$
with $\sigma'=\rho(a+ \sigma)$. If we choose
$\sigma> a\rho/(1-\rho)$ then $P$ becomes a strict cone contraction
of the regular cone $\calC_{\sigma,\RR}$. We also get the estimate
(to obtain an a priori estimate for the contraction
one may here try to optimize for the value of $\sigma$)~:
\begin{equation}
\label{Delta P estimate}
 \frac{\Delta_P}{2}
    \leq  \log \;\frac{\sigma+\rho(\sigma+a)}
          {\sigma-\rho(\sigma+a)} +
     D\; \rho\, (\sigma+a) .
\end{equation}

 By Theorem \ref{Thm Birkhoff},
   $P\in L(X_\RR)$ has a spectral gap
(see \cite{Rue68,FS79} and also \cite{Liv95}).\\
  
 Returning  to the complex operator, $M_g$, let us fix $y,y'\in \Omega$ and
 the corresponding pairing of pre-images $\calP(y,y')$ as described above.
 Let $\phi\in \calC_{\sigma',\RR}^*$ and
 write $\la \ell_{y,y'}, M_g \phi\ra = \sum_j \la \mu_j(g),\phi\ra$
with 
   \[ \la \mu_j(g),\phi\ra \equiv
        e^{g(x_j)} \phi(x_j) - 
            e^{-\sigma d(y,y')+g(x'_j)} \phi(x'_j), \ \ j\in J .\]
In order to compare with the real operator, we define complex numbers $w_j$,
$j\in J$,
through the relation
\[ \la \mu_j(g),\phi\ra  =
   e^{i\; \Im\; g(x_j)}
    \,
    w_j \, \la \mu_j(\Re\; g),\phi\ra .\]
 Equivalently (when the denominator is non-zero)~:
 \[
   e^{i\; \Im\; g(x_j)} \, w_j = 
        \frac{e^{ g(x_j)} \phi(x_j) - 
            e^{-\sigma d(y,y')+g(x'_j)} \phi(x'_j)}
        {e^{\Re\; g(x_j)} \phi(x_j) - 
            e^{-\sigma d(y,y') + \Re\; g(x'_j)} \phi(x'_j)}
	    .\]
We may apply Lemma \ref{exp estimate} below
with the bounds $\Re (z_1-z_2) \geq 
 \left( \sigma - \rho(\sigma+a)\right) d(y,y')$ and
 $|\Im (z_1-z_2)| \leq \rho b  \;d(y,y')$ to deduce that
\begin{equation}
 | \Arg \ w_j | \leq s_0 
        \equiv \frac {\rho b} { \sigma  - \rho(\sigma+a)}.
 \ \ \ \mbox{and} \ \ \ 
     1\  \leq \ {|w_{j}|^2}\leq {1+s_0^2}.
     \mylabel{w prefactor}
     \end{equation}

Given $i,j\in J$ we have~:
\[ \la \mu_j(g),\phi\ra \overline{\la \mu_i(g),\phi\ra} =
    \left( e^{i (\Im\; g(x_j) - \Im\; g(x_i))} w_j \overline{w_i} \right)\ \ 
    \la \mu_j(\Re \; g),\phi\ra \  {\la \mu_i(\Re\; g),\phi\ra} .\]
The two last factors are real and non-negative
(because $\sigma-\rho(\sigma+a)>0$)
and the complex pre-factor belongs to the set
\[ A = \{ r e^{i u}: 1 \leq r \leq 1+s_0^2,\ \ |u| \leq \theta+2 s_0\}.\]
Summing over all indices we therefore obtain
 \[
  \la \ell_{y,y'}, M_g \phi_1 \ra
    \la \ell_{w,w'}, \barM_g \phi_2 \ra = 
    Z
  \la \ell_{y,y'}, P \phi_1 \ra
    \la \ell_{w,w'}, P \phi_2 \ra ,\]
 in which   $Z$ is an average of numbers in $A$ whence belongs
 to $\Conv (A)$, the convex hull of $A$.
 
When $\theta+2 s_0 < \pi/4$ we conclude that 
the bounds in Assumption \ref{main assumption} are verified for the
constants
 $\alpha=\cos (\theta+2 s_0)$, $\gamma=(1+s_0^2) \sin (\theta+2 s_0)$
 and $\beta= 1+s_0^2$. The spectral gap condition in Theorem \ref{C domin sp gap}
 then reads as follows~:
 \begin{equation}
 \label{sp gap cond}
  (1+s_0^2) \; \tan(\theta+2 s_0) \; \cosh \frac{\Delta_P}{2} < 1 .
 \end{equation}

 Now, in order to get a more tractable and explicit formula we make  the 
 following (not optimal) choice for $\sigma$~:

 \[ \sigma = \frac{2a\rho}{1-\rho} + \frac{1}{\rho D}.\]

 Then $\sigma'=\rho(a+\sigma) \leq  \frac{1+\rho}{2} \sigma$ so that
 $(\sigma+\rho(a+\sigma))/(\sigma-\rho(a+\sigma)) \leq (3+\rho)/(1-\rho)$. 
 Using (\ref{Delta P estimate}) we obtain 
 \[ \cosh \frac{\Delta_P}{2} \leq e^{\Delta_P/2} 
     = \frac{3+\rho}{1-\rho}
       \exp\left( 1 +2 a \,D \rho \frac{1+\rho}{1-\rho} \right)
        .\]

One also checks that  \ 
$(\theta+2s_0) \frac{4}{1-\rho}<1$ \  implies  \ 
that $(1+s_0^2) \tan (\theta + 2 s_0) \frac{3+\rho}{1-\rho} < 1$ \
so we may replace (\ref{sp gap cond})
 by the stronger condition

 \[ (\theta + 2\, s_0) \; \exp \left( 
       1 +  \rho \,\frac{1+\rho}{1-\rho}\; D\, a \right) 
     \frac{4}{1-\rho} <1 .\]
Finally inserting $s_0=\rho^2\, D \,b /(1-\rho + \rho^2\, D\, a)$
we obtain the claimed condition which is thus sufficient 
for a spectral gap. \Halmos \\

\begin{Remark} In the literature, one often includes a
statement on Gibbs measures as well. If we let $\lambda h \otimes \mu$ denote
the leading spectral projection of $P=M_{\Re\; g}$, 
then positivity of $P$ implies that
the `state'
$\phi \in X_\RR \mapsto \nu(\phi) = \mu(\phi h)$ is uniformly bounded 
with respect to $|\phi|_0$. By continuity, $\nu$  extends to a linear functional
on $C^0(\Omega)$. If, in addition, we assume $\Omega$ compact, then by Riesz, this functional defines
a Borel probability measure $d\nu$ on $\Omega$. The measure is invariant
and strongly mixing for $f$.
 It is  known as a Gibbs measure for $f$ and 
the weight $g$.
 This part of the theorem, however,
needs the partial ordering induced by the cone of 
positive continuous functions and does not extend to a complex setup
(in general, it is even false there).
\end{Remark}

In the proof we made use of the following complex estimate~:
\begin{Lemma}
\mylabel{exp estimate}
  Let $z_1,z_2\in \CC$ be such that $\Re\, z_1 > \Re \,z_2$ and define
  $w\in\CC$ through
   \[ e^{i \; \Im \; z_1} w \equiv
               \frac
	           {e^{z_1} - e^{z_2} }
	           {e^{\Re\; z_1} - e^{\Re\; z_2} }.
		   \]
Then
 \[ |\Arg\; w | \leq \frac 
         {| \Im \; (z_1-z_2) | }
         { \Re \; (z_1-z_2)  } 
	 \ \ \mbox{and} \ \ 
	1 \leq |w^2| \leq
         1 + \left( \frac { \Im \; (z_1-z_2)  }
         { \Re \; (z_1-z_2) }  \right)^2 .
	 \]
\end{Lemma}
Proof: 
Writing $t=\Re\; (z_1-z_2) >0$ and $s=\Im\; (z_1-z_2)$ we have~:
\[ w = \frac{1-e^{-t-is}}{1-e^{-t}} .\]
Taking real and imaginary parts,
$\Re\; w = \frac{1 - e^{-t} \cos s}{1-e^{-t}}$
and 
$\Im \; w = \frac{e^{-t} \sin s}{1-e^{-t}}$,
we get
$|w|^2 = 1 + \frac{\sin^2 (s/2)}{\sinh^2(t/2)} \leq 1+(\frac{s}{t})^2$.
Also $|\frac{\partial}{\partial s} \, \log w| 
  = |\frac{1}{w} \frac{\partial w}{\partial s}|
  = \frac{e^{-t}}{ |1-e^{-t-is}|}
   \leq\frac{ e^{-t}}{1-e^{-t}} \leq \frac{1}{t}$ 
 so that $|\Arg\; w| \leq \frac{|s|}{t}$.\Halmos\\

\def\theequation{\Alph{section}.\arabic{equation}}

\appendix
\section{Projective space}
\mylabel{a proj metric}
Let  $X$ be a complex Banach space. 
Given non-zero elements $x,y\in X^*\equiv X-\{0\}$ 
we write $x\sim y$   iff\  $\CC x=\CC y$. Let $\pi: X^* \rr X^*/\sim$
denote the quotient map and write $[x]=\CC^*x$ for the equivalence
class of $x\in X^*$.
We equip the quotient space $\pi(X^*)$  with the following metric
\begin{equation}
d_{\pi(X^*)}([x],[y])
	= \dist_H (\CC x \cap S, \CC y \cap S)
      =\inf \left\{ \left\| 
        \frac{\mu x}{\|\mu x\|}\!-\!
        \frac{\nu y}{\|\nu y\|} \right\| : \mu,\nu\in \CC^*\right\},
	\ \ x,y\in X^*
  \label{proj metric}
\end{equation}
in which $\dist_H$ is the Hausdorff distance between non-empty sets
and $S=S(X)$ is the unit-sphere.

\begin{Lemma}
 \mylabel{XP lower bound}\mbox{}
  \begin{enumerate}
 \item  Let $\calC\subset X$ be a $\CC$-cone of $K$-bounded sectional
 aperture. Then for all $x,y\in \calC^*$~:
  \[ d_{\pi(X^*)}([x],[y])\leq 2 K d_{\calC}(x,y).\]
 \item  Let $x\in X^*$, $r>0$ and set $V=\pi^{-1} B_{\pi(X^*)}([x],r)$.
 Then for all $y\in V^*$
 \[ d_V(x,y) \leq \log \frac{r+d_{\pi(X^*)}([x],[y])}{r-d_{\pi(X^*)}([x],[y])}.\]
 \end{enumerate}
\end{Lemma}
Proof: 
Using the inequality
\begin{equation}
  \left\|\frac{x}{\|x\|} - \frac{y}{\|y\|}\right\| \leq 2\|x-y\| 
  \min \left\{\frac{1}{\|x\|}
  ,\frac{1}{\|y\|} \right\},\ \  x,y\in X^* 
  \mylabel{xy inequality}
\end{equation}
we obtain from Lemma \ref{bounded sectional aperture}~:
\[ 
  \left\|\frac{x/\la m,x\ra }{\|x/\la m,x\ra\|}
   - \frac{y/\la m,y\ra}{\|y/\la m,y\ra\|}\right\| \leq 
    2 \|m\|
  \left\|\frac{x}{\la m,x\ra}
   - \frac{y}{\la m,y\ra}\right\|
     \leq 2K d_{\calC}(x,y)\]
and the first conclusion follows.
For the second claim, normalize so that $d_{\pi(X^*)}([x],[y])=\|x-y\|<r$
and $\|x\|=\|y\|=1$. Let
 $u_\lambda=\frac{1+\lambda}{2} x +
        \frac{1-\lambda}{2} y =
         x+ \frac{1-\lambda}{2} (y-x)$.
By (\ref{xy inequality}), $\|\frac{u_\lambda}{\|u_\lambda\|} -x \|
 \leq \frac{|1-\lambda|}{2}\|y-x\|$ which remains smaller than $r$
 when $|1-\lambda| < \frac{2 r}{\|x-y\|} \equiv 2R \in (2,+\infty]$.
 Then
 $d_V(x,y) \leq d_{B_\CC(1,2R)}(-1,1) = d_{\DD} (0,\frac{1}{R})= \log \frac{R+1}{R-1}$
\Halmos\\          

Given any two points $x,y\in\calC^*$
we may follow Kobayashi \cite{Kob67,Kob70}
and define a projective pseudo-distance between $x$ and $y$ through~:
\[ \tilde{d}_\calC(x,y) = 
    \inf\{ \sum d_\calC(x_i,x_{i+1}): x_0=x,x_1,\ldots,x_n=y\in \calC^* \}.\]
Since $d_{\pi(X^*)}$ is a (projective) metric,
the previous Lemma implies that 
\begin{Theorem}
 \mylabel{2K Lipschitz}
 Suppose that $\calC$ is of $K$-bounded sectional aperture in $X$. Then
 the inclusion map,
  $(\calC^*, \tilde{d}_\calC) \rr (\calC^*, d_{\pi(X^*)})$ is $2K$-Lipschitz.
\end{Theorem}
In other words, this new distance
does not degenerate when taking the inf over finite chains,
so distinct complex
lines in $\calC$ have a non-zero $\tilde{d}_\calC$-distance.
This is conceptually very nice,
but, in our context, not particularly useful.
The reason is that even if $T\in L(X)$ maps $\calC^*$ into 
a subset of finite diameter in $\calC^*$ for the metric $\tilde{d}$ 
this does not seem to
imply a uniform contraction of $T$, i.e.\ no spectral gap.
We  leave a further study of this metric to the interested reader.\\

%\section{Open questions and future development}
%\begin{enumerate}
%\item Let $\CCC$ be the complexification of a real cone $\CCR$ and
%   consider  a complex cone-contraction $A: \CCC \rr \CCC$ having a unique
%   invariant complex line, $\CC h \subset \CCC$.
%    When $A$ is real then this invariant line is, in fact, simply a real
%    line times the complex field. When $A$ is `close' to real 
%   it may be of interest to know how `close' the invariant complex
%   line is to a real line (times $\CC$).
%   The interested reader  may  study this
%    e.g.\ through the contraction of the
%    following family of complexified cones
%   $\calC_{\CC,\theta} = ...$. The parameter $\theta$ indicates precisely how
%   far we are from the real case, with $\theta=0$ being the real case and
%   $\theta=\pi/2$ being the natural complexification studied in this paper.
% \item In our complex Kre\u{\i}n-Rutman Theorem we needed an auxiliary
%   assumption on the cone to assure strict positivity of the spectral 
%   radius of $A$. Is this really needed ?
% \item .....
%\end{enumerate}
%

\end{document}